\documentclass[10pt]{IEEEconf}

\usepackage{array}
\usepackage{latexsym}
\usepackage{amssymb}
\usepackage{amsmath}
\usepackage{amsfonts}
\usepackage{graphics}
\usepackage{layout}
\usepackage{epsf}
\usepackage{epsfig,color}

{\end{list}}

\newtheorem{Theorem}{Theorem}[section]
\newtheorem{Definition}[Theorem]{Definition}

\newtheorem{Proposition}[Theorem]{Proposition}

\begin{document}
\pdfminorversion 7

\title{\bf Networked Traffic State Estimation Involving Mixed Fixed-mobile Sensor Data Using Hamilton-Jacobi equations}


\author{Edward S. Canepa\thanks{E. Canepa is a Ph.D. Candidate, Department of Electrical Engineering, King Abdullah University of Science and Technology, Saudi Arabia (email: \texttt{edward.canepa@kaust.edu.sa}). Corresponding Author}, Christian G. Claudel\thanks{C. Claudel is an Assistant Professor,  Department of Civil, Architectural and Environmental Engineering, University of Texas Austin, USA (email:\texttt{christian.claudel@utexas.edu}).}}
\maketitle

\begin{abstract} 
Nowadays, traffic management has become a challenge for urban areas, which are covering larger geographic spaces and facing the generation of different kinds of traffic data. This article presents a robust traffic estimation framework for highways modeled by a system of Lighthill Whitham Richards equations that is able to assimilate different sensor data available. We first present an equivalent formulation of the problem using a Hamilton-Jacobi equation. Then, using a semi-analytic formula, we show that the model constraints resulting from the Hamilton-Jacobi equation are linear ones. We then pose the problem of estimating the traffic density given incomplete and inaccurate traffic data as a Mixed Integer Program. We then extend the density estimation framework to highway networks with any available data constraint and modeling junctions. Finally, we present a travel estimation application for a small network using real traffic measurements obtained obtained during \emph{Mobile Century} traffic experiment, and comparing the results with ground truth data. 
\end{abstract}

\section{Introduction}

Transportation research is currently at a tipping point; the emergence of new transformative technologies and systems, such as vehicle connectivity, automation, shared-mobility, and advanced sensing is rapidly changing the individual mobility and accessibility. This will fundamentally transform how transportation planning and operations should be conducted to enable smart and connected communities. The transport systems can be highly beneficiated and become safer, more efficient and reliable. Nowadays, dynamic routing and traffic-dependent navigation services are available for users. Such applications need to estimate the present traffic situation and that of the near future at a forecasting horizon based on data that are available in real-time. Traffic state estimation for a road network refers to estimate all the traffic variables (e.g. cars density, speed) of the network at an instant of time based of traffic measurements. This is, for a limited amount of traffic data the estimator obtains a complete view of the traffic scenario. This estimation requires the fusion or traffic data and traffic models, the latter are typically formulated as partial differential equations (PDEs). For this framework, we will use the the \emph{Lighthill-Whitham-Richards} (LWR) partial differential equation~\cite{LW55,Richards56} which is commonly used to model highway traffic; derivating the model constraints is a complex problem. Other estimation techniques such as Extended Kalman Filtering~\cite{AMSH04} (EKF), Ensemble Kalman Filtering~\cite{ARMXWB09} or Particle Filtering (PF) rely on approximations to determine the model constraints, either through linearization or sampling.

No approximation of the model is required by the framework presented on this article. An example of the usage of this framework is determining the ranges input flows (or any convex function of the boundary data) compatible with the traffic model and measurement data. The exact estimation technique presented in this paper is based on the \emph{Moskowitz} function~\cite{Moskowitz65, Newell93}; it is used here as an intermediate computational abstraction. The \emph{Moskowitz} function can be understood as the integral form of the density function, and solves an \emph{Hamilton-Jacobi} (HJ) PDE, whereas the density function itself solves the LWR PDE. An advantage of using the HJ PDE is that its solutions can be expressed semi analytically~\cite{TAC2}, which enables the derivation of the model constraints explicitly. 

\subsection{Contributions of the article}

The present article builds on \cite{CCITS12, NHM2013, Li2014525, TCNSLi, anderson2013optimization} which introduced a Mixed Integer Linear Programming framework for solving data assimilation and data reconciliaton problems, for specific objective functions. In the present article, we extend this framework to network traffic density estimation. The present article has the following contributions over existing work:

\begin{itemize}
\item The integration of internal traffic density data, or arbitrary travel time data (not necessarily defined as the travel time required to cross the entire physical domain), which was not considered in earlier articles.
\item The extension of the traffic state estimation framework defined earlier in \cite{NHM2013} to transportation networks, which require the proper modeling of junctions, and the integration of the entropy condition to junction flows.
\item The formulation of estimation problems that do not involve a minimum variance estimation, unlike classical estimation schemes derived from the Kalman Filter. Examples of non minimum variance estimation include compressed sensing ($L_1$ norm minimization), shown in Section \ref{s:highwaynet}.
\end{itemize}  


The outline of this article is the following. In Section~\ref{s:convexf} we define the solution to the LWR PDE and its equivalent formulation as a HJ PDE. In section~\ref{s:explsols}, we recall the analytical expressions of the solutions to HJ PDEs for the triangular flux functions investigated in this article, and show that the LWR PDE constraints correspond to convex constraints in the unknown initial, boundary and internal condition parameters. A first estimation example is shown in section~\ref{s:estprob}, using boundary and internal conditions from measurement data the unknown initial conditions are estimated. The framework is extended to Highway Networks in section~\ref{s:highwaynet}, where we also validated it using experimental traffic flow data (e.g. density, point velocity and travel time) collected during the \emph{Mobile Century} traffic experiment.

\section{Background}  \label{s:convexf}

\subsection{The Lighthill-Whitham-Richards traffic flow model}

For the remainder of the article, we will assume that the spatial domain representing the highway section is $[\xi,\chi]$, where $\xi$ and $\chi$ respectively represent the upstream and downstream boundaries of the domain.
Traffic flow on this section can be described by the \emph{density} function, denoted as $\rho(\cdot,\cdot)$. The density function represents an aggregated number of vehicles per space unit, and can is modeled by the \emph{Lighthill-Whitham-Richards} (LWR) PDE:

\begin{equation} \label{e:lwrpde} \frac{\partial \rho(t,x)}{\partial t} +\frac{\partial \psi(\rho(t,x))}{\partial x} \; = \; 0
 \end{equation}

The function $\psi(\cdot)$ is named \emph{flux function}. It depends on several empirical parameters (e.g. number of lanes, the drivers habits). Different models have been proposed for $\psi$, in
particular the triangular model defined below, this model is widely used in the literature~\cite{DaganzoCell94,Daganzo05,Daganzo06}.

\begin{equation} \label{e:triangh}
\psi(\rho)= \left\{ \begin{array}{ll}
v\rho & {\rm if} \;\;\rho \le \rho_c\\
w(\rho-\rho_m) & {\rm otherwise}
 \end{array} \right.
\end{equation}

In the remainder of this article, we assume that the flux function is triangular and given by~(\ref{e:triangh}). While other concave flux functions could be used and would also yield convex constraints, the instantiation of model constraints as linear inequalities requires piecewise linear flux functions, such as~(\ref{e:triangh}).

\subsection{Hamilton-Jacobi equation}

Equivalently, the state of traffic can be described by a scalar function ${\bf M}(\cdot,\cdot)$ of both time and space, known as~\emph{Moskowitz function}~\cite{Moskowitz65,Newell93}. The Moskowitz function is a macroscopic description of traffic flow, which appears naturally in the context of traffic. It can be interpreted as follows: let consecutive integer labels be assigned to vehicles entering the highway at location $x=\xi$. The Moskowitz function ${\bf M}(\cdot,\cdot)$ satisfies $\lfloor {\bf M}(t,x) \rfloor = n$ where $n$ is the label of the vehicle located in $x$ at time $t$~\cite{Daganzo05,Daganzo06}, and is assumed to be continuous.

The density function $\rho(\cdot,\cdot)$ is related~\cite{Newell93} to the spatial derivative of the Moskowitz function ${\bf M}(\cdot,\cdot)$ as follows:

\begin{equation}
\rho(t,x)= -\frac{\partial {\bf M}(t,x)}{\partial x}
\end{equation}

If the density function is to be modeled by the LWR PDE, the Moskowitz function satisfies an \emph{Hamilton-Jacobi} (HJ) PDE obtained~\cite{ABSP07,TAC1} by integration of the LWR PDE:

\begin{equation}\label{e:hjpde} \frac{\partial {\bf M}(t,x)}{\partial t}
    -\psi\left(-\frac{\partial {\bf M}(t,x)} {\partial x}\right)
 \; = \;  0 \end{equation}

Several classes of weak solutions to equation~(\ref{e:hjpde}) exist, such as viscosity solutions~\cite{CL83,BCD97} or Barron-Jensen/Frankowska (B-J/F) solutions~\cite{bj90hj,Frankowska93}. For the problem investigated in this article, these solutions are equivalent, and can be computed implicitly using a Lax-Hopf formula.

\subsection{Barron-Jensen/Frankowska solutions to Hamilton Jacobi equations}

In order to characterize the B-J/F solutions, we first need to define the Legendre-Fenchel transform of the Hamiltonian $\psi(\cdot)$ as follows.

\begin{Definition} \textbf{[Legendre-Fenchel transform]} For an upper semicontinuous Hamiltonian $\psi(\cdot)$, the Legendre-Fenchel transform $\varphi^*(\cdot)$ is given by:
\vspace{-0.05in} \begin{equation} \label{e:ctransform}
\varphi^*(u)\; :=  \begin{array}{ll} \displaystyle{\sup_{p \in \mbox{\rm Dom}(\psi )}[p \cdot u + \psi(p) ]} \\
\end{array}
\end{equation}
\end{Definition}

Solving the HJ PDE~(\ref{e:hjpde}) requires the definition of~\emph{value conditions}, which encode the traditional concepts of initial, boundary and internal conditions.

\begin{Definition} \label{def:valcond} \textbf{[Value condition]} A value condition ${\bf c}(\cdot,\cdot)$ is a lower semicontinuous function defined on a subset of $[0,t_{\max}] \times [\xi,\chi]$.
\end{Definition}

In the remainder of the article, a value condition can encode an initial condition, an upstream boundary condition or a downstream boundary condition. Each of these functions is defined on a subset of $\mathbb{R}_+ \times [\xi,\chi]$.

For each value condition ${\bf c}(\cdot,\cdot)$, we define the partial solution~\cite{MCB10TRB} to the HJ PDE~(\ref{e:hjpde}) using the Lax-Hopf formula~\cite{ABSP07,TAC1}.

\vspace{0.1in} \begin{Proposition} \label{p:LaxHopf} \textbf{[Lax-Hopf formula]} Let $\psi(\cdot)$ be a concave Hamiltonian, and let $\varphi^*(\cdot)$ be its Legendre-Fenchel transform~(\ref{e:ctransform}). Let ${\bf c}(\cdot,\cdot)$ be a lower semicontinuous value condition, as in Definition~\ref{def:valcond}. The B-J/F solution ${\bf M}_{\bf c}(\cdot,\cdot)$ to~(\ref{e:hjpde}) associated with ${\bf c}(\cdot,\cdot)$ can be algebraically represented~\cite{ABSP07,TAC1} by:

\begin{equation} \label{e:LaxHopf}
\footnotesize {\bf M}_{\bf c}(t,x) =  \displaystyle{\inf_{(u,T) \in \mbox{\rm
Dom}(\varphi^{*}) \times \mathbb{R}_+  }} \left( {\bf c} (t-T, x+ T u) + T \varphi^{*}(u )\right)
\end{equation}
\end{Proposition}

Equation~(\ref{e:LaxHopf}) implies the existence of a B-J/F solution ${\bf M}_{\bf c}(\cdot,\cdot)$ for any value condition function ${\bf c}(\cdot,\cdot)$. However, the solution itself may be incompatible with the value condition that we imposed on it, \emph{i.e.} we do not necessarily have $\forall (t,x) \in {\rm Dom}({\bf c}), {\bf M}_{\bf c}(t,x)={\bf c}(t,x)$.

\subsection{Properties of the solutions to scalar Hamilton-Jacobi equations}

The structure of the Lax-Hopf formula~(\ref{e:LaxHopf}), implies the following important property, known as~\emph{inf-morphism} property. The inf-morphism property can be formally derived through capture basins, such as in~\cite{ABSP07}.

\begin{Proposition} \textbf{[Inf-morphism property]} Let the value condition ${\bf c}(\cdot,\cdot)$ be minimum of a finite number of lower semicontinuous functions:
\begin{equation} \label{e:definfmorp}
\forall (t,x) \in [0,t_{\max}] \times [\xi,\chi], \;\; {\bf c}(t,x):= \displaystyle{\min_{j \in J} {\bf c}_j(t,x) }
\end{equation}

The solution~${\bf M}_{\bf c}(\cdot,\cdot)$ associated with the above value condition can be decomposed~\cite{ABSP07,TAC1,TAC2} as:

\begin{equation} \label{e:definfmorp2}
\forall (t,x) \in [0,t_{\max}] \times [\xi,\chi], \;\; {\bf M}_{\bf c}(t,x)= \displaystyle{\min_{j\in J} {\bf M}_{{\bf c}_j}(t,x)}
\end{equation}
\end{Proposition}

In the present work, the value conditions ${\bf c}_j(\cdot,\cdot)$ are not known exactly, either because of measurement uncertainty (case of the upstream and downstream boundary condition) or because of the lack of measurements (case of the initial condition). However, even if the real values of ${\bf c}_j(\cdot,\cdot)$ are not known exactly, they cannot be arbitrary as they have to apply in the strong sense (see~\cite{SB06} for a mathematical definition) to be compatible with the LWR model. In the remainder of this article, we define the model constraints as the set of constraints that applies on the value conditions ${\bf c}_j(\cdot,\cdot)$ to ensure that all value conditions apply in the strong sense.

The inf-morphism property is critical for the derivation of the LWR PDE model constraints, allowing us to instantiate these constraints as inequalities.

\begin{Proposition} \textbf{[Model compatibility constraints for block value conditions]} Let ${\bf c}(\cdot,\cdot)=\displaystyle{\min_{j \in J} {\bf c}_j(\cdot,\cdot)}$ be given, and let ${\bf M}_{\bf c}(\cdot,\cdot)$ be defined as in~(\ref{e:LaxHopf}). The value condition ${\bf c}(\cdot,\cdot)$ satisfies $\forall (t,x) \in {\rm Dom}({\bf c}), {\bf M}_{\bf c}(t,x)={\bf c}(t,x)$ if and only if the following inequality constraints are satisfied:
\begin{equation} \label{e:compconst}
{\bf M}_{{\bf c}_j}(t,x) \ge {\bf c}_i(t,x) \;\;\forall (t,x) \in {\rm Dom}({\bf c}_i),\;\; \forall (i,j) \in J^2
\end{equation}
\end{Proposition}

The proof of this proposition is available in~\cite{SIAM11CB}. Note that equation~(\ref{e:compconst}) represents an important improvement, as the model constraints are now semi-explicit. In order to solve the problem completely, we still need to evaluate the functions  ${\bf M}_{{\bf c}_j}(\cdot,\cdot)$ explicitly. These explicit solutions were derived in~\cite{TAC2} for affine initial, boundary and internal conditions blocks which are presented in the next section.

\section{Explicit solutions to piecewise affine initial, internal and boundary conditions} \label{s:explsols}

As mentioned before, different types of value conditions can be incorporated into the estimation problem. In the present article  we will handle initial, internal, upstream and downstream conditions. These value conditions are typically measured (with some error) using fixed sensors, such as inductive loop detectors, magnetometers or GPS devices.

\subsection{Definition of affine initial, upstream/downstream boundary and internal density conditions}

The formal definition of initial, upstream/downstream boundary and internal conditions associated with the HJ PDE~(\ref{e:hjpde}) is the subject of the following definition.

\begin{Definition} \textbf{[Affine initial, upstream/downstream boundary and internal conditions]} Let us define $\mathbb{K}=\{0,\dots,k_{\max}\}$, $\mathbb{N}=\{0,\dots,n_{\max}\}$, $\mathbb{M}=\{0,\dots,m_{\max}\}$ and $\mathbb{U}=\{0,\dots,u_{\max}\}$. For all $k \in \mathbb{K}$, $n \in \mathbb{N}$, $m \in \mathbb{M}$ and $u \in \mathbb{U}$, we define the following functions, respectively called initial, upstream, downstream internal flow and internal density conditions:

\begin{equation} \label{e:defvalcondini}
\footnotesize
\begin{array}{ll}
M_k(t,x) \hspace{-0.02in}= \hspace{-0.02in} \left\{ \hspace{-0.05in}
	\begin{array}{ll}
	-\sum_{i=0}^{k-1}{\rho_{\rm ini}(i)}X \\ -\rho_{\rm ini}(k)(x-kX) & {\rm if }\;\; t=0 \\ & {\rm and } \;\;x \in[kX,(k+1)X]\\
	+\infty & {\rm otherwise}
	\end{array}\right.
\end{array}
\end{equation}
\begin{equation} \label{e:defvalcondu}
\footnotesize
\begin{array}{ll}
\gamma_n(t,x) \hspace{-0.02in}= \hspace{-0.02in} \left\{ \hspace{-0.05in}
	\begin{array}{ll}
	\sum_{i=0}^{n-1}{q_{\rm in}(i)}T \\ +q_{\rm in}(n)(t-nT) & {\rm if }\;\; x=\xi \\ & {\rm and } \;\;t \in[nT,(n+1)T]\\
	+\infty & {\rm otherwise}
	\end{array}\right.
\end{array}
\end{equation}
\begin{equation}\label{e:defvalcondd} \footnotesize
\begin{array}{ll}
\beta_n(t,x)\hspace{-0.02in} =\hspace{-0.02in} \left\{ \hspace{-0.05in}
	\begin{array}{ll}
	\sum_{i=0}^{n-1}{q_{\rm out}(i)}T\\+q_{\rm out}(n)(t-nT)\\ -\sum_{k=0}^{k_{max}} \rho(k)X  & {\rm if }\;\; x=\chi \\ & {\rm and }\; t \in[nT,(n+1)T]\\
	+\infty & {\rm otherwise}
	\end{array}\right.
\end{array}
\end{equation}
\begin{equation}\label{e:defvalcondi} \footnotesize
\begin{array}{ll}
\mu_m(t,x) \hspace{-0.02in}= \hspace{-0.02in} \left\{ \hspace{-0.05in}
	\begin{array}{ll}
	L(m)+r(m)(t-t_{\rm min}(m)) &{\rm if }\;\; x=x_{\min}(m)\\ & +v^{\rm meas}(m)(t-t_{\rm min}(m)) \\ & {\rm and }\; t \in[t_{\rm min}(m),t_{\rm max}(m)]\\
	+\infty & {\rm otherwise}
	\end{array}\right.
\end{array}
\end{equation}
\begin{equation}\label{e:defvalcondu} \footnotesize
\begin{array}{ll}
\Upsilon_u(t,x) \hspace{-0.02in}= \hspace{-0.02in} \left\{ \hspace{-0.05in}
	\begin{array}{ll}
	L(u)-\rho(u)(x-x_{\rm min_{\rho}}(u)) &{\rm if }\;\; x \in[x_{\min_{\rho}}(u),x_{\max_{\rho}}(u)]\\ & {\rm and }\; t =t_{\rho}(u)\\
	+\infty & {\rm otherwise}
	\end{array}\right.
\end{array}
\end{equation}

where $v^{\rm meas}(m)=\frac{x_{\rm max}(m)-x_{\rm min}(m)}{t_{\rm max}(m)-t_{\rm min}(m)}$

\end{Definition}

\begin{figure}[h]
\begin{center} \begin{tabular}{c}
{\mbox{\includegraphics[width = 0.85\linewidth]{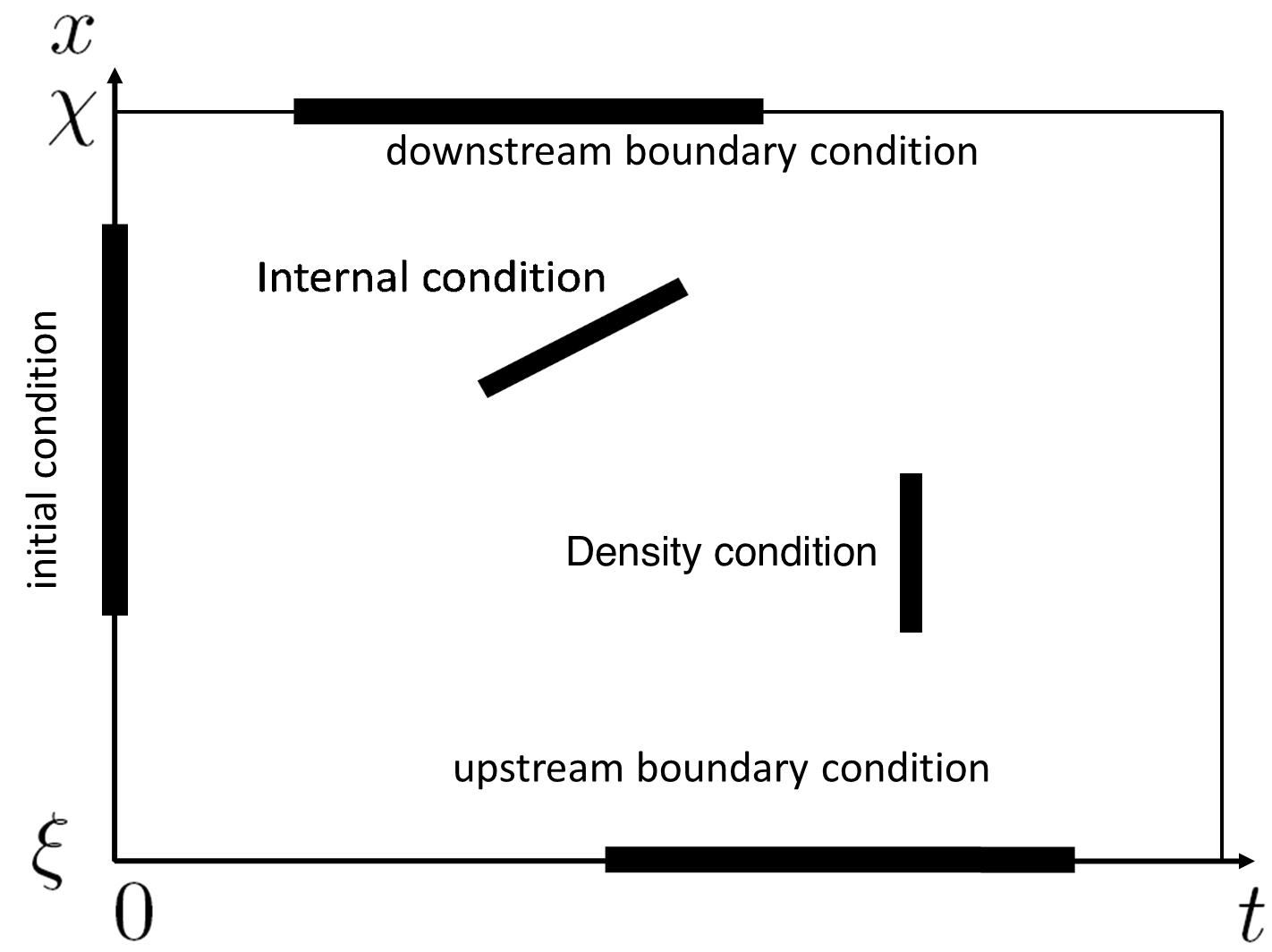}}}
\end{tabular}
\end{center}
\vspace{-0.1in} \caption{\textbf{Domains of the initial, upstream/downstream boundary and internal conditions.} The block upstream and downstream boundary conditions respectively denoted by $\gamma_n(\cdot,\cdot)$ and $\beta_n(\cdot,\cdot)$ are defined on line segments corresponding to the upstream and downstream boundaries of the physical domain. In contrast, the block initial conditions $M_k(\cdot,\cdot)$ are defined on line segments corresponding to the initial time. Note that the actual problem involves block initial conditions covering the entire physical domain $[\xi,\chi]$, and block boundary conditions covering the temporal domain $[0,t_{\max}]$. The block of internal flow and density conditions ($\mu_m(\cdot,\cdot)$ and $\Upsilon_u(\cdot,\cdot)$ respectively) are defined and any position inside the computational domain. All these functions are unknown (or only partially known).\label{f:illdom}}
\end{figure}

In the above definition, internal density conditions (\ref{e:defvalcondu}) are specific to model density sensors that are inside our computational domain. Regarding flow sensors also located inside the computational domain, they can be thought as an internal flow condition (\ref{e:defvalcondi}) associated with a zero velocity ($v^{\rm meas}(m)=0$). Note that the affine initial, upstream/downstream boundary and internal conditions defined above for the HJ PDE~(\ref{e:hjpde}) are equivalent to constant initial, upstream/downstream boundary and internal conditions for the LWR PDE~(\ref{e:lwrpde}).
The domains of definitions of these functions are illustrated in Figure~\ref{f:illdom}.

\subsection{Analytical solutions to affine initial, upstream/downstream boundary and internal conditions}

Given the affine initial, upstream/downstream boundary and internal conditions defined above, the corresponding solutions ${\bf M}_{M_k}(\cdot,\cdot)$, ${\bf M}_{\gamma_n}(\cdot,\cdot)$, ${\bf M}_{\beta_n}(\cdot,\cdot)$ and and ${\bf M}_{\Upsilon_u}(\cdot,\cdot)$ are given~\cite{SIAM11CB,MCB10TRB} by the following formulas:

\begin{equation} \label{e:solutions}
\footnotesize
\begin{array}{ll}
{\bf M}_{M_k}(t,x)=
\begin{cases}
	+\infty &{\rm if } \;\; x \le kX+wt \\ &{\rm or} \;\; x \ge (k+1)X+vt\\
	-\sum_{i=0}^{k-1}{\rho_{\rm ini}(i)}X \\+\rho_{\rm ini}(k)(tv+kX-x)& {\rm if } \;\;kX+tv\leq x\\ \hfill  {\rm and} & (k+1)X +tv\geq x\\ \hfill {\rm and} & \rho_{\rm ini}(k) \le \rho_c\\
    -\sum_{i=0}^{k-1}{\rho_{\rm ini}(i)}X\\ +\rho_c(tv+kX-x)& {\rm if } \;\;kX+tv\geq x\\  \hfill {\rm and} & kX +tw \leq x\\ \hfill {\rm and} & \rho_{\rm ini}(k) \le \rho_c\\
	-\sum_{i=0}^{k-1}{\rho_{\rm ini}(i)}X \\+\rho_{\rm ini}(k)(tw+kX-x)\\-\rho_{m}tw& {\rm if } \;\;kX+tw\leq x\\ \hfill  {\rm and} & (k+1)X +tw\geq x\\ \hfill {\rm and} & \rho_{\rm ini}(k) \ge \rho_c\\
    -\sum_{i=0}^{k}{\rho_{\rm ini}(i)}X\\ \rho_c(tw+(k+1)X-x)\\-\rho_{m}tw& {\rm if } \;\;(k+1)X+tv\geq x\\  \hfill {\rm and} & (k+1)X +tw \leq x\\ \hfill {\rm and} & \rho_{\rm ini}(k) \ge \rho_c\\
	\end{cases}\\
\end{array}
\end{equation}

\begin{equation} \label{e:solutions2}
\footnotesize
\begin{array}{ll}
{\bf M}_{\gamma_n}(t,x)=
\begin{cases}
	+\infty &{\rm if } \;\; t\leq nT+\frac{x-\xi}{v}\\
	\sum_{i=0}^{n-1}{q_{\rm in}(i)T}\\+q_{\rm in}(n)(t-\frac{x-\xi}{v}-nT)& {\rm if } \;\;nT+\frac{x-\xi}{v}\leq t\\ &  {\rm and}\;\; t\leq (n+1)T\\ &+\frac{x-\xi}{v}\\
	\sum_{i=0}^{n}{q_{\rm in}(i)T}\\+\rho_c v(t-(n+1)T-\frac{x-\xi}{v})&{\rm otherwise}\\
	\end{cases}\\
\end{array}
\end{equation}
\begin{equation} \label{e:solutions3}
\footnotesize
\begin{array}{ll}
{\bf M}_{\beta_n}(t,x)=
\begin{cases}
	+\infty &{\rm if }\;\; t\leq nT\\ & +\frac{x-\chi}{w}\\
	-\sum_{k=0}^{k_{max}} \rho_{\rm ini}(k)X+\sum_{i=0}^{n-1}{q_{\rm out}(i)T}\\+q_{\rm out}(n)(t-\frac{x-\chi}{w}-nT) \\-\rho_m(x-\chi)&{\rm if }\;\; nT\\ & +\frac{x-\chi}{w}\leq t\\ \hfill {\rm and} & t\leq (n+1)T\\ & +\frac{x-\chi}{w}\\
	-\sum_{k=0}^{k_{max}} \rho(k)X+\sum_{i=0}^{n}{q_{\rm out}(i)T}\\+\rho_c v(t-(n+1)T-\frac{x-\chi}{v})& {\rm otherwise}\\
	\end{cases}\\
\end{array}
\end{equation}

\begin{equation}\label{e:solutions4}
\footnotesize
\begin{array}{ll}
{\bf M}_{\mu_m}(t,x)\hspace{-0.02in} =\hspace{-0.02in}
\begin{cases}
L_m+\\r_m\left(t-\frac{x-x_{\rm min}(m)-v^{\rm meas}(m)(t-t_{\rm min}(m))}{v-v^{\rm meas}(m)}-t_{\rm min}(m)\right)\\
\mbox{    if }x\geq x_{\rm min}(m)+v^{\rm meas}(m)(t-t_{\rm min}(m))\\
\mbox{    and }x\geq x_{\rm max}(m)+v(t-t_{\rm max}(m))\\
\mbox{    and }x\leq x_{\rm min}(m)+v(t-t_{\rm min}(m))\\
L_m+\\r_m\left(t-\frac{x-x_{\rm min}(m)-v^{\rm meas}(m)(t-t_{\rm min}(m))}{w-v^{\rm meas}(m)}-t_{\rm min}(m)\right)\\ +k_c(v-w)\frac{x-x_{\rm min}(m)-v^{\rm meas}(m)(t-t_{\rm min}(m))}{w-v^{\rm meas}(m)}\\
\mbox{    if }x\leq x_{\rm min}(m)+v^{\rm meas}(m)(t-t_{\rm min}(m))\\
\mbox{    and }x\leq x_{\rm max}(m)+w(t-t_{\rm max}(m))\\
\mbox{    and }x\geq x_{\rm min}(m)+w(t-t_{\rm min}(m))\\
L_m+r_m\left(t_{\rm max}(m)-t_{\rm min}(m)\right)+\\ \left(t-t_{\rm max}(m)\right)k_c\left(v-\frac{x-x_{\rm max}(m)}{t-t_{\rm max}(m)}\right)\\
\mbox{    if }x\leq x_{\rm max}(m)+v(t-t_{\rm max}(m))\\
\mbox{    and }x\geq x_{\rm max}(m)+w(t-t_{\rm max}(m))\\
+\infty\mbox{  otherwise}
\end{cases}
\end{array}
\end{equation}

\begin{equation} \label{e:solutions5}
\footnotesize
\begin{array}{ll}
{\bf M}_{\Upsilon_u}(t,x)=
\begin{cases}
	+\infty &{\rm if } \;\; x \le x_{\rm min_{\rho}}(u)\\ & + w(t-t_{\rho}(u)) \\ &{\rm or} \;\; x \ge x_{\rm max_{\rho}}(u)\\ & + v(t-t_{\rho}(u))\\ &{\rm or}
\;\; t \le t_{\rho}(u)\\	
	L(u) \\+\rho(u)(v(t-t_{\rm \rho}(u))+x_{\rm min_{\rho}}(u)-x)& {\rm if } \;x_{\rm min_{\rho}}(u)\\ & + v(t-t_{\rho}(u))\leq x\\ \hfill  {\rm and} & x_{\rm max_{\rho}}(u)\\ & +v(t-t_{\rho}(u))\geq x\\ \hfill {\rm and} & \rho(u) \le \rho_c\\
    	L(u) \\+\rho_c(v(t-t_{\rho}(u))+x_{\rm min_{\rho}}(u)-x)& {\rm if } \;\;x_{\rm min_{\rho}}(u)\\ &+w(t-t_{\rho}(u))\leq x\\ \hfill  {\rm and} & x_{\rm min_{\rho}}(u)\\ &+v(t-t_{\rho}(u))\geq x\\ \hfill {\rm and} & \rho(u) \le \rho_c\\
 	L(u) \\+\rho(u)(w(t-t_{\rho}(u))+x_{\rm min_{\rho}}(u)-x)\\ -\rho_{\rm m}w(t-t_{\rho}(u))& {\rm if } \;\;x_{\rm min_{\rho}}(u)\\ &+w(t-t_{\rho}(u))\leq x\\ \hfill  {\rm and} & x_{\rm max_{\rho}}(u)\\ &+w(t-t_{\rho}(u))\geq x\\ \hfill {\rm and} & \rho(u) \geq \rho_c\\
    L(u) \\+\rho_{\rm c}(w(t-t_{\rho}(u))+x_{\rm max_{\rho}}-x)\\ -\rho_{\rm m}w(t-t_{\rho}(u))& {\rm if } \;\;x_{\rm max_{\rho}}(u)\\ &+w(t-t_{\rho}(u))\leq x\\ \hfill  {\rm and} & x_{\rm max_{\rho}}(u)\\ &+v(t-t_{\rho}(u))\geq x\\ \hfill {\rm and} & \rho(u) \geq \rho_c\\
	\end{cases}\\
\end{array}
\end{equation}

\subsection{Properties of the affine initial, upstream/downstream boundary and internal conditions}

In this section, we show that the LWR model constraints~(\ref{e:compconst}) are convex in the variable $\Big(\rho(1),\rho(2),\dots,\rho(k_{\max}),q_{\rm in}(1),\dots,q_{\rm in}(n_{\max}),q_{\rm out}(1),$ $\dots,q_{\rm out}(n_{\max}) \Big)$.

\begin{Proposition} \label{p:conc1} \textbf{[Linearity property of the initial, upstream/downstream boundary and internal conditions]} Let us fix $(t,x) \in \mathbb{R}_+ \times [\xi,\chi]$. The initial, upstream and downstream boundary condition functions $M_{k}(t,x)$, $\gamma_n(t,x)$ and $\beta_n(t,x)$ are linear functions of the coefficients $(\rho(1),\dots,\rho(k_{\max}),q_{\rm in }(1),\dots,q_{\rm in}(n_{\max}),q_{\rm out}(1),$ $ \dots,q_{\rm out}(n_{\max}))$.
\end{Proposition}

The proof of this proposition is straightforward and follows directly from equations~(\ref{e:defvalcondini}), (\ref{e:defvalcondu}) and~(\ref{e:defvalcondd}).

\begin{Proposition} \label{p:conc2} \textbf{[Concavity property of the solution associated with the initial condition]} Let us fix $(t,x) \in \mathbb{R}_+ \times [\xi,\chi]$. The solution ${\bf M}_{{ M}_{k}}(t,x)$ associated with the initial condition~(\ref{e:defvalcondini}) is a concave function of the coefficients $\rho(\cdot)$.
\end{Proposition}

{\bf Proof} --- \hspace{ 2 mm} The Lax-Hopf formula~(\ref{e:LaxHopf}) associated with the solution ${\bf M}_{M_k}(\cdot,\cdot)$ can be written~\cite{TAC1,TAC2} as:

\begin{equation}\label{e:laxHopfic}\footnotesize {\bf M}_{M_k}(t,x)\; =\; \begin{array}{c} \displaystyle{\inf_{u \in \mbox{\rm Dom}(\varphi^{*})\;\;\text{s. t.} \;\;(x+tu) \in [kX,(k+1)X]}}
 \Big( -\sum_{i=0}^{k-1}{\rho(i)}X  \\  -\rho(k)(x-kX)+ t\varphi^{*}(u)\Big) \end{array} \end{equation}

Let us fix $u \times {\rm Dom}(\varphi^*)$. The function $f$ defined as $f(\rho(1),\dots,\rho(k_{\max}))= -\sum_{i=0}^{k-1}{\rho(i)}X-\rho(k)(x-kX)+ t\varphi^{*}(u)$ is concave (indeed, affine). Hence, the function ${\bf M}_{M_k}(t,x)$ is a concave function of $(\rho(1),\dots,\rho(k_{\max}))$, since it is the infimum of concave functions~\cite{Boyd,rockafellar1970ca}. \hfill $\;\; \blacksquare$

\begin{Proposition} \label{p:conc3} \textbf{[Concavity property of the solutions associated with upstream and downstream boundary conditions]} Let us fix $(t,x) \in \mathbb{R}_+ \times [\xi,\chi]$. The solutions ${\bf M}_{\gamma_n}(t,x)$ and ${\bf M}_{\beta_n}(t,x)$ respectively associated with the upstream and downstream boundary conditions~(\ref{e:defvalcondu}) and~(\ref{e:defvalcondd}) are concave functions of the coefficients $\rho(\cdot)$,  $q_{\rm in}(\cdot)$ and $q_{\rm out}(\cdot)$.
\end{Proposition}

{\bf Proof} --- \hspace{ 2 mm} The Lax-Hopf formula~(\ref{e:LaxHopf}) associated with the solution ${\bf M}_{\gamma_n}(\cdot,\cdot)$ can be written~\cite{TAC1,TAC2} as:

\begin{equation} \label{e:lxibcaffp} \begin{array}{l}
{\bf M}_{\gamma_n}(t,x)= \displaystyle{\inf_{ s \in \mathbb{R}_+ \cap [t-(n+1)T,t-nT]}}
\sum_{i=0}^{n-1}{q_{\rm in}(i)}T\\ +q_{\rm in}(n)(t-s)  +s \varphi^{*}(\frac{\xi-x}{s})
\end{array} \end{equation}

Let us fix $s \in \mathbb{R}_+ \cap [t-(n+1)T,t-nT]$. The function $d$ defined as $d(q_{\rm in}(1),\dots,q_{\rm in}(n_{\max}))= \sum_{i=0}^{n-1}{q_{\rm in}(i)}T +q_{\rm in}(n)(t-s)  +s \varphi^{*}(\frac{\xi-x}{s})$ is concave (indeed, affine). Hence, the solution ${\bf M}_{\gamma_n}(t,x)$ is a concave function of $(q_{\rm in}(1),\dots,q_{\rm in}(n_{\max}))$, since it is the infimum of concave functions~\cite{Boyd,rockafellar1970ca}.
The same property applies for ${\bf M}_{\beta_n}(t,x)$, which is a concave function of $(\rho(1),\dots,\rho({k_{\max}}),q_{\rm out}(1),\dots,q_{\rm out}(n_{\max}))$ \hfill $\;\; \blacksquare$

Propositions~\ref{p:conc1},~\ref{p:conc2} and~\ref{p:conc3} thus imply the following convexity property:

\begin{Proposition} \label{p:modelconst} \textbf{[Convexity property of model constraints]} The model constraints~(\ref{e:compconst}) are convex functions of $\big(\rho(1),\rho(2),\dots,\rho(k_{\max}),q_{\rm in}(1),\dots,q_{\rm in}(n_{\max}),q_{\rm out}(1),$ $\dots,q_{\rm out}(n_{\max}) \big)$.
\end{Proposition}

{\bf Proof} --- \hspace{ 2 mm} The set of inequality constraints~(\ref{e:compconst}) can be written as:

\begin{equation} \label{e:ec1}
\begin{array}{c} {\bf M}_{{\bf c}_i}(t,x) \ge {\bf c}_j(t,x), \;\; \forall (t,x) \in {\rm Dom}({\bf c}_j)\\ \;\; \forall j \in I \text{ such that } (t,x) \in {\rm Dom}({\bf c}_j),\;\; \forall i \in I \end{array}
\end{equation}

Note that Proposition~\ref{p:conc1} implies that the term ${\bf c}_j(t,x)$ in~(\ref{e:ec1}) is a linear function (labeled $l_{j,t,x}(\cdot)$) of $\big(\rho(1),\rho(2),\dots,\rho(k_{\max}),q_{\rm in}(1),\dots,q_{\rm in}(n_{\max}),q_{\rm out}(1),$ $\dots,q_{\rm out}(n_{\max}) \big)$. In addition, by Propositions~\ref{p:conc2} and~\ref{p:conc3}, the term ${\bf M}_{{\bf c}_i}(t,x)$ is a concave function (labeled $c_{i,t,x}(\cdot)$) of $\big(\rho(1),\rho(2),\dots,\rho(k_{\max}),q_{\rm in}(1),\dots,q_{\rm in}(n_{\max}),q_{\rm out}(1),$ $\dots,q_{\rm out}(n_{\max}) \big)$. Hence, the equality~(\ref{e:ec1}) can be written as:

\begin{equation} \label{e:ec2}
\begin{array}{l}
-c_{i,t,x}\Big(\rho(1),\rho(2),\dots,\rho(k_{\max}),q_{\rm in}(1),\dots,q_{\rm in}(n_{\max}),\\q_{\rm out}(1),\dots,q_{\rm out}(n_{\max}) \Big)+l_{j,t,x}\Big(\rho(1),\rho(2),\dots,\rho(k_{\max}),\\q_{\rm in}(1),\dots,q_{\rm in}(n_{\max}),q_{\rm out}(1),$ $\dots,q_{\rm out}(n_{\max}) \Big) \le 0, \\ \forall j \in I, \forall (t,x) \in {\rm Dom}({\bf c}_j),\;\; \forall i \in I
\end{array}
\end{equation}

This last inequality is a convex inequality~\cite{Boyd} in $\big(\rho(1),\rho(2),\dots,\rho(k_{\max}),q_{\rm in}(1),\dots,q_{\rm in}(n_{\max}),q_{\rm out}(1),$ $\dots,q_{\rm out}(n_{\max}) \big)$, that is, an inequality of the form $f(\cdot) \le 0$ where $f(\cdot)$ is a convex function. \hfill $\;\; \blacksquare$

The above property is very important, and can be thought of as follows. Consider the vector space ${\cal V}$ of all parameters of the initial, upstream and downstream boundary conditions. Each point of this vector space corresponds to a known value condition (encompassing initial, upstream and downstream boundary conditions). However, the solution to the LWR PDE~(\ref{e:lwrpde}) associated with this arbitrary value condition will satisfy the value condition itself on its boundaries if and only if the model constraints~(\ref{e:compconst}) are satisfied. Proposition~\ref{p:modelconst} essentially states that the set of value conditions compatible with the LWR PDE model is convex.

\section{Formulation of the density estimation problem as a Mixed Integer Linear Program} \label{s:estprob}

\subsubsection{Decision variable}

As outlined in Proposition~\ref{p:modelconst}, the variable $\big(\rho(1),\rho(2),\dots,\rho(k_{\max}),q_{\rm in}(1),\dots,q_{\rm in}(n_{\max}),q_{\rm out}(1),$ $\dots,q_{\rm out}(n_{\max}) \big)$ plays an important role in our estimation problem, and will be defined as the decision variable of our optimization framework.

\begin{Definition} \textbf{[Decision variable]} \label{def:defdecv} Let us consider a finite set of,
initial, upstream and downstream boundary conditions be defined as in~(\ref{e:defvalcondini})~(\ref{e:defvalcondu}) and~(\ref{e:defvalcondd}). The decision variable $v$ associated with this finite set of value conditions is defined by:

\begin{equation} \label{e:defdecv}
\begin{array}{l}
v:=\Big(\rho(1),\rho(2),\dots,\rho(k_{\max}),q_{\rm in}(1),\dots,q_{\rm in}(n_{\max}),\\q_{\rm out}(1),\dots,q_{\rm out}(n_{\max}) \Big)
\end{array}
\end{equation}
\end{Definition}

We denote by ${\cal V}$ the vector space of the decision variables $v$ defined by equation~(\ref{e:defdecv}).

\subsubsection{Model and data constraints}

Let $\overline{v}$ denote the value of the decision variable associated with the true state of the system (which is not known in practice, and can only be estimated). Because of model and data constraints, $\overline{v}$ must satisfy the set of constraints outlined in Propositions~\ref{p:modelconstdef} and~\ref{p:dataconst} below.

\begin{Proposition} \label{p:modelconstdef} \textbf{[Model constraints]} The model constraints~(\ref{e:compconst}) can be expressed as the following finite set of convex inequality constraints:

\end{Proposition}

\begin{equation}\label{e:ineqconsts1a}
\footnotesize
\begin{cases}
{\bf M}_{M_k}(t_{0},x_{p})\geq M_p(t_{0},x_{p}) \\ \hfill\forall (k,p)\in\mathbb{K}^2\quad (i)\\
{\bf M}_{M_k}(pT,\chi)\geq \beta_p(pT,\chi) \\ \hfill \forall (k,p)\in\mathbb{K}\times\mathbb{N}\quad (ii)(a)\\
{\bf M}_{M_k}(t_{0}+\frac{\chi-x_{k+1}}{v},\chi)\geq \beta_p(t_{0}+\frac{\chi-x_{k+1}}{v},\chi) \\ \quad \quad \hfill\forall (k,p)\in\mathbb{K}\times\mathbb{N} \;\; {\rm s.}\;\; {\rm t.} \;\;  t_{0} + \frac{\chi-x_{k+1}}{v} \in [pT, (p+1)T]  \quad(ii) (b)\\
{\bf M}_{M_k}(pT,\xi)\geq \gamma_p(pT,\xi) \\ \hfill \forall (k,p)\in\mathbb{K}\times\mathbb{N}\quad (iii)(a)\\
{\bf M}_{M_k}(t_{0}+\frac{\xi-x_{k}}{w},\xi)\geq \gamma_p(t_{0}+\frac{\xi-x_{k}}{w},\xi) \\ \hfill\forall (k,p)\in\mathbb{K}\times\mathbb{N} \;\; {\rm s.}\;\; {\rm t.} \;\;  t_{0} + \frac{\xi-x_{k}}{w} \in [pT, (p+1)T]  \quad(iii) (b)\\
\end{cases}
\end{equation}

\begin{equation}\label{e:ineqconsts1b}
\footnotesize
\begin{cases}
{\bf M}_{M_k}(t_{min}(m),x_{min}(m))\geq \mu_m(t_{min}(m),x_{min}(m)) \\ \hfill \forall k\in\mathbb{K}, \forall m  \in\mathbb{M} \quad (iv)(a)   \\
{\bf M}_{M_k}(t_{max}(m),x_{max}(m))\geq \mu_m(t_{max}(m),x_{max}(m)) \\ \hfill \forall k\in\mathbb{K}, \forall m \in\mathbb{M} \quad (iv)(b)  \\
{\bf M}_{M_k}(t_{1}(m,k),x_{1}(m,k))\geq \mu_m(t_{1}(m,k),x_{1}(m,k)) \\ \quad \quad \forall k\in\mathbb{K}, \forall m \in\mathbb{M} \; {\rm s.}\; {\rm t.}\;  \forall t_{1}(m,k) \in  [t_{min}(m), t_{max}(m)]  \quad (iv)(c)   \\
{\bf M}_{M_k}(t_{2}(m,k),x_{2}(m,k))\geq \mu_m(t_{2}(m,k),x_{2}(m,k))  \\  \quad \quad \forall k\in\mathbb{K}, \forall m \in\mathbb{M} \; {\rm s.}\; {\rm t.}\;  \forall t_{2}(m,k) \in  [t_{min}(m), t_{max}(m)]  \quad (iv)(d)   \\
{\bf M}_{M_k}(t_{3}(m,k),x_{3}(m,k))\geq \mu_m(t_{3}(m,k),x_{3}(m,k))  \\  \quad \quad \forall k\in\mathbb{K}, \forall m \in\mathbb{M} \; {\rm s.}\; {\rm t.}\;  \forall t_{3}(m,k) \in  [t_{min}(m), t_{max}(m)]  \quad (iv)(e)   \\
{\bf M}_{M_k}(t_{4}(m,k),x_{4}(m,k))\geq \mu_m(t_{4}(m,k),x_{4}(m,k))  \\  \quad \quad \forall k\in\mathbb{K}, \forall m \in\mathbb{M} \; {\rm s.}\; {\rm t.}\;  \forall t_{4}(m,k) \in  [t_{min}(m), t_{max}(m)]  \quad (iv)(f)   \\

\end{cases}
\end{equation}

\begin{equation}\label{e:ineqconsts1c}
\footnotesize
\begin{cases}
{\bf M}_{M_k}(t_{\rho}(u),x_{\min_{\rho}}(u))\geq \Upsilon_u(t_{\rho}(u),x_{\min_{\rho}}(u)) \\ \hfill \forall k\in\mathbb{K}, \forall u  \in\mathbb{U} \quad (v)(a)   \\
{\bf M}_{M_k}(t_{\rho}(u),x_{\max_{\rho}}(u))\geq \Upsilon_u(t_{\rho}(u),x_{\max_{\rho}}(u)) \\ \hfill \forall k\in\mathbb{K}, \forall u  \in\mathbb{U} \quad (v)(b)   \\
{\bf M}_{M_k}(t_{\rho}(u),x_{5}(u,k))\geq \Upsilon_u(t_{\rho}(u),x_{5}(u,k)) \\ \quad \quad \forall k\in\mathbb{K}, \forall u \in\mathbb{U} \; {\rm s.}\; {\rm t.}\;  x_{5}(u,k) \in  [x_{min_{\rho}}(u), x_{max_{\rho}}(u)]  \quad (v)(c)   \\
{\bf M}_{M_k}(t_{\rho}(u),x_{6}(u,k))\geq \Upsilon_u(t_{\rho}(u),x_{6}(u,k)) \\ \quad \quad \forall k\in\mathbb{K}, \forall u \in\mathbb{U} \; {\rm s.}\; {\rm t.}\;  x_{6}(u,k) \in  [x_{min_{\rho}}(u), x_{max_{\rho}}(u)]  \quad (v)(d)   \\
{\bf M}_{M_k}(t_{\rho}(u),x_{7}(u,k))\geq \Upsilon_u(t_{\rho}(u),x_{7}(u,k)) \\ \quad \quad \forall k\in\mathbb{K}, \forall u \in\mathbb{U} \; {\rm s.}\; {\rm t.}\;  x_{7}(u,k) \in  [x_{min_{\rho}}(u), x_{max_{\rho}}(u)]  \quad (v)(e)\\  
{\bf M}_{M_k}(t_{\rho}(u),x_{8}(u,k))\geq \Upsilon_u(t_{\rho}(u),x_{8}(u,k)) \\ \quad \quad \forall k\in\mathbb{K}, \forall u \in\mathbb{U} \; {\rm s.}\; {\rm t.}\;  x_{8}(u,k) \in  [x_{min_{\rho}}(u), x_{max_{\rho}}(u)]  \quad (v)(f)  \\

\end{cases}
\end{equation}

\begin{equation} \label{e:ineqconsts2a}
\footnotesize
\begin{cases}
{\bf M}_{\gamma_n}(pT,\xi)\geq\gamma_p(pT,\xi)  & \forall (n,p)\in\mathbb{N}^2\hfill (vi)\\
{\bf M}_{\gamma_n}(pT,\chi)\geq\beta_p(pT,\chi) &  \forall (n,p)\in\mathbb{N}^2\hfill (vii) (a)\\
{\bf M}_{\gamma_n}(nT+\frac{\chi-\xi}{v},\chi)\geq \beta_p(nT+\frac{\chi-\xi}{v},\chi) & \forall (n,p)\in\mathbb{N}^2 \;\; {\rm s.}\;\; {\rm t.} \;\; nT+\\ & \frac{\chi-\xi}{v} \in [pT,(p+1)T] \\ & \hfill  (vii) (b)\\
\end{cases}
\end{equation}
\begin{equation} \label{e:ineqconsts2b}
\footnotesize
\begin{cases}
{\bf M}_{\gamma_n}(t_{\min}(m),x_{\min}(m))\geq \mu_m(t_{\min}(m),x_{\min}(m)) \\ \hfill \forall n\in\mathbb{N},\forall m\in\mathbb{M}\quad (viii) (a)\\
{\bf M}_{\gamma_n}(t_{\max}(m),x_{\max}(m))\geq \mu_m(t_{\max}(m),x_{\max}(m)) \\ \hfill \forall n\in\mathbb{N},\forall m\in\mathbb{M}\quad (viii) (b)\\
{\bf M}_{\gamma_n}(t_9(m,n),x_9(m,n))\geq\mu_m(t_9(m,n),x_9(m,n)) \\ \quad \quad \forall n\in\mathbb{N},\forall m\in\mathbb{M}\;\; {\rm s.} \;\; {\rm t.}\;\; t_9(m,n) \in [t_{\min}(m);t_{\max}(m)] \quad (viii) (c)\\
\end{cases}
\end{equation}

\begin{equation} \label{e:ineqconsts2c}
\footnotesize
\begin{cases}
{\bf M}_{\gamma_n}(t_{\rho}(u),x_{\min_{\rho}}(u))\geq \Upsilon_u(t_{\rho}(u),x_{\min_{\rho}}(u)) \\ \hfill \forall n\in\mathbb{N},\forall u\in\mathbb{U}\quad (ix) (a)\\
{\bf M}_{\gamma_n}(t_{\rho}(u),x_{\max_{\rho}}(u))\geq \Upsilon_u(t_{\rho}(u),x_{\max_{\rho}}(u)) \\ \hfill \forall n\in\mathbb{N},\forall u\in\mathbb{U}\quad (ix) (b)\\
{\bf M}_{\gamma_n}(t_{\rho}(u),x_{10}(u,n))\geq\Upsilon_u(t_{\rho}(u),x_{10}(u,n)) \\ \quad \quad \forall n\in\mathbb{N},\forall u\in\mathbb{U}\;\; {\rm s.} \;\; {\rm t.}\;\; x_{10}(u,n) \in [x_{\min_{\rho}}(u);x_{\max_{\rho}}(u)] \quad (ix) (c)\\
\end{cases}
\end{equation}

\begin{equation} \label{e:ineqconsts3a}
\footnotesize
\begin{cases}
{\bf M}_{\beta_n}(pT,\xi)\geq\gamma_p(pT,\xi) \\ \hfill \forall (n,p)\in\mathbb{N}^2 \quad (x)(a)\\
{\bf M}_{\beta_n}(nT+\frac{\xi-\chi}{w},\xi)\geq\gamma_p(nT+\frac{\xi-\chi}{w},\xi)\\ \quad \quad \forall (n,p)\in\mathbb{N}^2 \;\; {\rm s.} \;\; {\rm t.} \;\; nT+\frac{\xi-\chi}{w} \in [pT,(p+1)T] \quad (x)(b)\\
{\bf M}_{\beta_n}(pT,\chi)\geq\beta_p(pT,\chi) \\ \hfill \forall (n,p)\in\mathbb{N}^2 \quad (xi)\\
\end{cases}
\end{equation}
\begin{equation}\label{e:ineqconsts3b}
\footnotesize
\begin{cases}
{\bf M}_{\beta_n}(t_{\min}(m),x_{\min}(m))\geq\mu_m(t_{\min}(m),x_{\min}(m)) \\ \hfill \forall n\in\mathbb{N},\forall m\in\mathbb{M}\quad (xii) (a)\\
{\bf M}_{\beta_n}(t_{\max}(m),x_{\max}(m))\geq\mu_m(t_{\max}(m),x_{\max}(m)) \\ \hfill \forall n\in\mathbb{N},\forall m\in\mathbb{M}\quad (xii) (b)\\
{\bf M}_{\beta_n}(t_{11}(m,n),x_{11}(m,n))\geq\mu_m(t_{11}(m,n),x_{11}(m,n)) \\ \quad \quad \forall n\in\mathbb{N},\forall m\in\mathbb{M}\;\; {\rm s.} \;\; {\rm t.}  \;\;t_{11}(m,n) \in [t_{\min}(m);t_{\max}(m)]\;\;(xii) (c)\\
\end{cases}
\end{equation}

\begin{equation}\label{e:ineqconsts3c}
\footnotesize
\begin{cases}
{\bf M}_{\beta_n}(t_{\rho}(u),x_{\min_{\rho}}(u))\geq\Upsilon_u(t_{\rho}(u),x_{\min_{\rho}}(u)) \\ \hfill \forall n\in\mathbb{N},\forall u\in\mathbb{U}\quad (xiii)(a)\\
{\bf M}_{\beta_n}(t_{\rho}(u),x_{\max_{\rho}}(u))\geq\Upsilon_u(t_{\rho}(u),x_{\max_{\rho}}(u)) \\ \hfill \forall n\in\mathbb{N},\forall u\in\mathbb{U}\quad (xiii)(b)\\
{\bf M}_{\beta_n}(t_{\rho}(u),x_{12}(u,n))\geq\Upsilon_u(t_{\rho}(u),x_{12}(u,n)) \\ \quad \quad \forall n\in\mathbb{N},\forall m\in\mathbb{M}\;\; {\rm s.} \;\; {\rm t.}  \;\;x_{12}(m,n) \in [x_{\min_{\rho}}(u);x_{\max_{\rho}}(u)]\;\;(xiii)(c)\\
\end{cases}
\end{equation}

\begin{equation}\label{e:ineqconsts4a}
\footnotesize
\begin{cases}
{\bf M}_{\mu_m}(pT,\xi)\geq\gamma_p(pT,\xi) & \forall (m,p)\in\mathbb{M}\times\mathbb{N} \hfill (xiv) (a)\\
{\bf M}_{\mu_m}(t_{13}(m),\xi)\geq\gamma_p(t_{13}(m),\xi) &\forall (m,p)\in\mathbb{M}\times\mathbb{N} \; {\rm s.}\; {\rm t.} \\ & \hfill t_{13}(m) \in [pT,(p+1)T]\quad (xiv) (b)\\
{\bf M}_{\mu_m}(t_{14}(m),\xi)\geq\gamma_p(t_{14}(m),\xi) &\forall (m,p)\in\mathbb{M}\times\mathbb{N} \; {\rm s.}\; {\rm t.} \\ & \hfill t_{14}(m) \in [pT,(p+1)T]\quad (xiv) (c)\\
\end{cases}
\end{equation}

\begin{equation}\label{e:ineqconsts4b}
\footnotesize
\begin{cases}
{\bf M}_{\mu_m}(pT,\chi)\geq\beta_p(pT,\chi) & \hfill\forall (m,p)\in\mathbb{M}\times\mathbb{N}\quad (xv)(a)\\
{\bf M}_{\mu_m}(t_{15}(m),\chi)\geq\beta_p(t_{15}(m),\chi) &\forall (m,p)\in\mathbb{M}\times\mathbb{N}\; {\rm s.}\; {\rm t.}\\ & \hfill t_9(m) \in [pT,(p+1)T]\quad (xv)(b)\\
{\bf M}_{\mu_m}(t_{16}(m),\chi)\geq\beta_p(t_{16}(m),\chi) &\forall (m,p)\in\mathbb{M}\times\mathbb{N}\; {\rm s.}\; {\rm t.}\\ & \hfill t_{10}(m) \in [pT,(p+1)T]\quad (xv)(c)\\
\end{cases}
\end{equation}

\begin{equation}\label{e:ineqconsts4c}
\footnotesize
\begin{cases}
{\bf M}_{\mu_m}(t_{\min}(p),x_{\min}(p))\geq\mu_p(t_{\min}(p),x_{\min}(p)) \\ \hfill  \forall (m,p)\in\mathbb{M}^2\quad (xvi)(a)\\
{\bf M}_{\mu_m}(t_{\max}(p),x_{\max}(p))\geq\mu_p(t_{\max}(p),x_{\max}(p)) \\ \hfill \forall (m,p)\in\mathbb{M}^2\quad (xvi)(b)\\
{\bf M}_{\mu_m}(t_{17}(m,p),x_{17}(m,p))\geq\mu_p(t_{17}(m,p),x_{17}(m,p)) \\ \hfill \forall (m,p)\in\mathbb{M}^2 \;{\rm s.}\; {\rm t.}\; t_{17}(m,p) \in [t_{\min}(p),t_{\max}(p)] \quad (xvi)(c)\\
{\bf M}_{\mu_m}(t_{18}(m,p),x_{18}(m,p))\geq\mu_p(t_{18}(m,p),x_{18}(m,p)) \\ \hfill \forall (m,p)\in\mathbb{M}^2 \;{\rm s.}\; {\rm t.}\; t_{18}(m,p) \in [t_{\min}(p),t_{\max}(p)] \quad (xvi)(d)\\
{\bf M}_{\mu_m}(t_{19}(m,p),x_{19}(m,p))\geq\mu_p(t_{19}(m,p),x_{19}(m,p)) \\ \hfill \forall (m,p)\in\mathbb{M}^2 \;{\rm s.}\; {\rm t.}\; t_{19}(m,p) \in [t_{\min}(p),t_{\max}(p)] \quad (xvi)(e)\\
{\bf M}_{\mu_m}(t_{20}(m,p),x_{20}(m,p))\geq\mu_p(t_{20}(m,p),x_{20}(m,p)) \\ \hfill \forall (m,p)\in\mathbb{M}^2 \;{\rm s.}\; {\rm t.}\; t_{20}(m,p) \in [t_{\min}(p),t_{\max}(p)] \quad (xvi)(f)\\
{\bf M}_{\mu_m}(t_{21}(m,p),x_{21}(m,p))\geq\mu_p(t_{21}(m,p),x_{21}(m,p)) \\ \quad \quad \quad  \forall (m,p)\in\mathbb{M}^2 \;{\rm s.}\; {\rm t.}\; t_{21}(m,p) \in [t_{\min}(p),t_{\max}(p)] \quad (xvi)(g)\\
\end{cases}
\end{equation}

\begin{equation}\label{e:ineqconsts4d}
\footnotesize
\begin{cases}
{\bf M}_{\mu_m}(t_{\rho}(u),x_{\min_{\rho}}(u))\geq\Upsilon_u(t_{\rho}(u),x_{\min_{\rho}}(u)) \\ \hfill  \forall m\in\mathbb{M}, \forall u\in\mathbb{U}\quad (xvii)(a)\\
{\bf M}_{\mu_m}(t_{\rho}(u),x_{\max_{\rho}}(u))\geq\Upsilon_u(t_{\rho}(u),x_{\max_{\rho}}(u)) \\ \hfill  \forall m\in\mathbb{M}, \forall u\in\mathbb{U}\quad (xvii)(b)\\
{\bf M}_{\mu_m}(t_{\rho}(u),x_{22}(m,u))\geq\Upsilon_u(t_{\rho}(u),x_{22}(m,u)) \\ \hfill \forall m\in\mathbb{M}, \forall u\in\mathbb{U} \;{\rm s.}\; {\rm t.}\; x_{22}(m,u) \in [x_{\min_{\rho}}(u),x_{\max_{\rho}}(u)] \quad (xvii)(c)\\
{\bf M}_{\mu_m}(t_{\rho}(u),x_{23}(m,u))\geq\Upsilon_u(t_{\rho}(u),x_{23}(m,u)) \\ \hfill \forall m\in\mathbb{M}, \forall u\in\mathbb{U} \;{\rm s.}\; {\rm t.}\; x_{23}(m,u) \in [x_{\min_{\rho}}(u),x_{\max_{\rho}}(u)] \quad (xvii)(d)\\
{\bf M}_{\mu_m}(t_{\rho}(u),x_{24}(m,u))\geq\Upsilon_u(t_{\rho}(u),x_{24}(m,u)) \\ \hfill \forall m\in\mathbb{M}, \forall u\in\mathbb{U} \;{\rm s.}\; {\rm t.}\; x_{24}(m,u) \in [x_{\min_{\rho}}(u),x_{\max_{\rho}}(u)] \quad (xvii)(e)\\
{\bf M}_{\mu_m}(t_{\rho}(u),x_{25}(m,u))\geq\Upsilon_u(t_{\rho}(u),x_{25}(m,u)) \\ \hfill \forall m\in\mathbb{M}, \forall u\in\mathbb{U} \;{\rm s.}\; {\rm t.}\; x_{25}(m,u) \in [x_{\min_{\rho}}(u),x_{\max_{\rho}}(u)] \quad (xvii)(f)\\
{\bf M}_{\mu_m}(t_{\rho}(u),x_{26}(m,u))\geq\Upsilon_u(t_{\rho}(u),x_{26}(m,u)) \\ \hfill \forall m\in\mathbb{M}, \forall u\in\mathbb{U} \;{\rm s.}\; {\rm t.}\; x_{26}(m,u) \in [x_{\min_{\rho}}(u),x_{\max_{\rho}}(u)] \quad (xvii)(g)\\

\end{cases}
\end{equation}

\begin{equation}\label{e:ineqconsts5a}
\footnotesize
\begin{cases}
{\bf M}_{\Upsilon_u}(pT,\chi)\geq \beta_p(pT,\chi) \\ \hfill\forall (u,p)\in\mathbb{U}\times\mathbb{N} \quad(xviii)(a)\\
{\bf M}_{\Upsilon_u}(t_{\rho}(u)+\frac{\chi-x_{max_{\rho}}(u)}{v},\chi)\geq \beta_p(t_{\rho}(u)+\frac{\chi-x_{max_{\rho}}(u)}{v},\chi) \\ \forall (u,p)\in\mathbb{U}\times\mathbb{N} \;\; {\rm s.}\;\; {\rm t.} \;\;  t_{\rho}(u) + \frac{\chi-x_{max_{\rho}}(u)}{v} \in [pT,(p+1)T]  \quad(xviii)(b)\\
{\bf M}_{\Upsilon_u}(pT,\xi)\geq \gamma_p(pT,\xi) \\ \hfill \forall (u,p)\in\mathbb{U}\times\mathbb{N}\quad (xviii)(c)\\
{\bf M}_{\Upsilon_u}(t_{\rho}(u)+\frac{\xi-x_{min_{\rho}}(u)}{w},\xi)\geq \gamma_p(t_{\rho}(u)+\frac{\xi-x_{min_{\rho}}(u)}{w},\xi) \\ \hfill \forall (u,p)\in\mathbb{K}\times\mathbb{N} \;\; {\rm s.}\;\; {\rm t.} \;\;  t_{\rho}(u)+\frac{\xi-x_{min_{\rho}}(u)}{w} \in [pT,(p+1)T]  \quad(xviii)(d)\\
\end{cases}
\end{equation}

\begin{equation}\label{e:ineqconsts5b}
\footnotesize
\begin{cases}
{\bf M}_{\Upsilon_u}(t_{\min}(p),x_{\min}(p))\geq\mu_p(t_{\min}(p),x_{\min}(p)) \\ \hfill  \forall (u,p)\in\mathbb{U}\times\mathbb{M}\quad (xix)(a)\\
{\bf M}_{\Upsilon_u}(t_{\max}(p),x_{\max}(p))\geq\mu_p(t_{\max}(p),x_{\max}(p)) \\ \hfill \forall (m,p)\in\mathbb{U}\times\mathbb{M}\quad (xix)(b)\\
{\bf M}_{\Upsilon_u}(t_{27}(u,p),x_{27}(u,p))\geq\mu_p(t_{27}(u,p),x_{27}(u,p)) \\ \hfill \forall (u,p)\in\mathbb{U}\times\mathbb{M} \;{\rm s.}\; {\rm t.}\; t_{27}(u,p) \in [t_{\min}(p),t_{\max}(p)] \quad (xix)(c)\\
{\bf M}_{\Upsilon_u}(t_{28}(u,p),x_{28}(u,p))\geq\mu_p(t_{28}(u,p),x_{28}(u,p)) \\ \hfill \forall (u,p)\in\mathbb{U}\times\mathbb{M} \;{\rm s.}\; {\rm t.}\; t_{28}(u,p) \in [t_{\min}(p),t_{\max}(p)] \quad (xix) (d)\\
{\bf M}_{\Upsilon_u}(t_{29}(u,p),x_{29}(u,p))\geq\mu_p(t_{29}(u,p),x_{29}(u,p)) \\ \hfill \forall (u,p)\in\mathbb{U}\times\mathbb{M} \;{\rm s.}\; {\rm t.}\; t_{29}(u,p) \in [t_{\min}(p),t_{\max}(p)] \quad (xix) (e)\\
{\bf M}_{\Upsilon_u}(t_{30}(u,p),x_{30}(u,p))\geq\mu_p(t_{30}(u,p),x_{30}(u,p)) \\ \hfill \forall (u,p)\in\mathbb{U}\times\mathbb{M} \;{\rm s.}\; {\rm t.}\; t_{30}(u,p) \in [t_{\min}(p),t_{\max}(p)] \quad (xix) (f)\\
{\bf M}_{\Upsilon_u}(t_{\rho}(u),x_{26}(p,u))\geq\mu_p(t_{\rho}(u),x_{26}(p,u)) \\ \hfill \forall (u,p)\in\mathbb{U}\times\mathbb{M} \;{\rm s.}\; {\rm t.}\; x_{26}(p,u) \in [x_{\min_{\rho}}(u),t_{\max_{\rho}}(u)] \quad (xix) (g)\\
\end{cases}
\end{equation}

\begin{equation}\label{e:ineqconsts5c}
\footnotesize
\begin{cases}
{\bf M}_{\Upsilon_u}(t_{\rho}(p),x_{\min_{\rho}}(p))\geq\Upsilon_p(t_{\rho}(p),x_{\min_{\rho}}(p)) \\ \hfill  \forall (u,p)\in\mathbb{U}^2\quad (xx) (a)\\
{\bf M}_{\Upsilon_u}(t_{\rho}(p),x_{\max_{\rho}}(p))\geq\Upsilon_p(t_{\rho}(p),x_{\max_{\rho}}(p)) \\ \hfill  \forall (u,p)\in\mathbb{U}^2\quad (xx) (b)\\
{\bf M}_{\Upsilon_u}(t_{\rho}(p),x_{31}(u,p))\geq\mu_p(t_{\rho}(p),x_{31}(u,p)) \\ \hfill \forall (u,p)\in\mathbb{U}^2 \;{\rm s.}\; {\rm t.}\; x_{31}(u,p) \in [x_{\min_{\rho}}(p),t_{\max_{\rho}}(p)] \quad (xx) (c)\\
{\bf M}_{\Upsilon_u}(t_{\rho}(p),x_{32}(u,p))\geq\mu_p(t_{\rho}(p),x_{32}(u,p)) \\ \hfill \forall (u,p)\in\mathbb{U}^2 \;{\rm s.}\; {\rm t.}\; x_{32}(u,p) \in [x_{\min_{\rho}}(p),t_{\max_{\rho}}(p)] \quad (xx) (d)\\
{\bf M}_{\Upsilon_u}(t_{\rho}(p),x_{33}(u,p))\geq\mu_p(t_{\rho}(p),x_{33}(u,p)) \\ \hfill \forall (u,p)\in\mathbb{U}^2 \;{\rm s.}\; {\rm t.}\; x_{33}(u,p) \in [x_{\min_{\rho}}(p),t_{\max_{\rho}}(p)] \quad (xx) (e)\\
{\bf M}_{\Upsilon_u}(t_{\rho}(p),x_{34}(u,p))\geq\mu_p(t_{\rho}(p),x_{34}(u,p)) \\ \hfill \forall (u,p)\in\mathbb{U}^2 \;{\rm s.}\; {\rm t.}\; x_{34}(u,p) \in [x_{\min_{\rho}}(p),t_{\max_{\rho}}(p)] \quad (xx) (f)\\

\end{cases}
\end{equation}

where the coefficients $t_1(m,k)$, $x_1(m,k)$, $t_2(m,k)$, $x_2(m,k)$, $t_3(m,k)$, $x_3(m,k)$, $t_4(m,k)$, $x_4(m,k)$, $x_5(u,k)$, $x_6(u,k)$, $x_7(u,k)$, $x_8(u,k)$, $t_9(m,n)$, $x_9(m,n)$, $x_{10}(u,n)$, $t_{11}(m,n)$, $x_{11}(m,n)$, $x_{12}(u,n)$, $t_{13}(m)$, $t_{14}(m)$, $t_{15}(m)$, $t_{16}(m)$, $t_{17}(m,p)$, $x_{17}(m,p)$, $t_{18}(m,p)$, $x_{18}(m,p)$, $t_{19}(m,p)$, $x_{19}(m,p)$, $t_{20}(m,p)$, $x_{20}(m,p)$, $t_{21}(m,p)$, $x_{21}(m,p)$, $x_{22}(m,u)$, $x_{23}(m,u)$, $x_{24}(m,u)$, $x_{25}(m,u)$, $x_{26}(m,u)$, $t_{27}(u,p)$, $x_{27}(u,p)$, $t_{28}(u,p)$, $x_{28}(u,p)$, $t_{29}(u,p)$, $x_{29}(u,p)$, $t_{30}(u,p)$, $x_{30}(u,p)$, $x_{31}(u,p)$, $x_{32}(u,p)$, $x_{33}(u,p)$ and $x_{34}(u,p)$ are given by equations~(\ref{e:constants1}),~(\ref{e:constants2}),~(\ref{e:constants3}) and~(\ref{e:constants4}) below:

\begin{equation} \label{e:constants1}
\footnotesize
\left\{ \begin{array}{ll}
t_1(m,k)=\frac{x_{\min}(m)-x_{\rm k + 1} - v^{\rm meas}(m)t_{\min}(m)}{v-v^{\rm meas}(m)}\\
x_1(m,k)=v^{\rm meas}(m) \Big(\frac{x_{\min}(m)-x_{\rm k+1}-v_{\rm meas}(m)t_{\min}(m) }{v-v_{\rm meas}(m)} \\ - t_{\min}(m)\Big)  + x_{\min}(m)  \\
t_2(m,k)=\frac{x_{\min}(m)-x_{\rm k } - v^{\rm meas}(m)t_{\min}(m)}{w-v^{\rm meas}(m)}\\
x_2(m,k)=v^{\rm meas}(m) \Big(\frac{x_{\min}(m)-x_{\rm k}-v_{\rm meas}(m)t_{\min}(m) }{w-v_{\rm meas}(m)} \\ - t_{\min}(m)\Big)   + x_{\min}(m)\\
t_3(m,k)=\frac{x_{\min}(m)-x_{\rm k } - v^{\rm meas}(m)t_{\min}(m)}{v-v^{\rm meas}(m)}\\
x_3(m,k)=v^{\rm meas}(m) \Big(\frac{x_{\min}(m)-x_{\rm k}-v_{\rm meas}(m)t_{\min}(m) }{v-v_{\rm meas}(m)} \\ - t_{\min}(m)\Big)  + x_{\min}(m)  \\
t_4(m,k)=\frac{x_{\min}(m)-x_{\rm k +1} - v^{\rm meas}(m)t_{\min}(m)}{w-v^{\rm meas}(m)}\\
x_4(m,k)=v^{\rm meas}(m) \Big(\frac{x_{\min}(m)-x_{\rm k}-v_{\rm meas}(m)t_{\min}(m) }{w-v_{\rm meas}(m)} \\ - t_{\min}(m)\Big)   + x_{\min}(m)\\
x_5(u,k)=x_{k+1}+v(t_{\rho}(u)-t_{0})\\
x_6(u,k)=x_{k+1}+w(t_{\rho}(u)-t_{0})\\
x_7(u,k)=x_{k}+v(t_{\rho}(u)-t_{0})\\
x_8(u,k)=x_{k}+w(t_{\rho}(u)-t_{0})\\
\end{array} \right.
\end{equation}

\begin{equation} \label{e:constants2}
\footnotesize
\left\{ \begin{array}{ll}
t_9(m,n)=\frac{nTv-v^{\rm meas}(m)t_{\rm min}(m)+x_{min}(m)-\xi}{v-v^{\rm meas}(m)}\\
x_9(m,n)=x_{\min}(m)+ \\v^{\rm meas}(m)\left(\frac{nTv-v^{\rm meas}(m)t_{\min}(m)+x_{\min}(m)-\xi}{v-v^{\rm meas}(m)}-t_{\min}(m)\right) \\
x_{10}(u,n) = \xi + v(t_{\rho}(u)-nT)\\
t_{11}(m,n)=\frac{nTw-v^{\rm meas}(m)t_{\min}(m)+x_{\min}(m)-\chi}{w-v^{\rm meas}(m)}\\
x_{11}(m,n)=x_{\min}(m)+ \\v^{\rm meas}(m)\left(\frac{nTw-v^{\rm meas}(m)t_{\min}(m)+x_{\min}(m)-\chi}{w-v^{\rm meas}(m)}-t_{\min}(m)\right)\\
x_{12}(u,n) = \chi + w(t_{\rho}(u)-nT)\\
t_{13}(m)=\frac{\xi-x_{\min}(m)+wt_{\min}(m)}{w}\\
t_{14}(m)=\frac{\xi-x_{\max}(m)+wt_{\max}(m)}{w}\\
t_{15}(m)=\frac{\chi-x_{\min}(m)+vt_{\min}(m)}{v}\\
t_{16}(m)=\frac{\chi-x_{\max}(m)+vt_{\max}(m)}{v}\\
\end{array} \right.
\end{equation}

and

\begin{equation}\label{e:constants3}
\footnotesize
\left\{ \begin{array}{ll}
t_{17}(m,p)=\frac{x_{\min}(m)-x_{\min}(p)+v^{\rm meas}(p)t_{\min}(p)-v^{\rm meas}(m)t_{\min}(m)}{v^{\rm meas}(p)-v^{\rm meas}(m)}\\
x_{17}(m,p)=x_{\min}(p)+v^{\rm meas}(p)\Big(-t_{\rm min}(p) +  \\  \frac{x_{\min}(m)-x_{\min}(p)+v^{\rm meas}(p)t_{\min}(p)-v^{\rm meas}(m)t_{\min}(m)}{v^{\rm meas}(p)-v^{\rm meas}(m)}\Big)\\
t_{18}(m,p)=\frac{x_{\max}(m)-x_{\min}(p)+v^{\rm meas}(p)t_{\min}(p)-vt_{\max}(m)}{v^{\rm meas}(p)-v}\\
x_{18}(m,p)=x_{\min}(p)+v^{\rm meas}(p)\Big(-t_{\rm min}(p)+ \\ \frac{x_{\max}(m)-x_{\min}(p)+v^{\rm meas}(p)t_{\min}(p)-vt_{\max}(m)}{v^{\rm meas}(p)-v}  \Big)\\
t_{19}(m,p)=\frac{x_{\min}(m)-x_{\min}(p)+v^{\rm meas}(p)t_{\min}(p)-vt_{\min}(m)}{v^{\rm meas}(p)-v}\\
x_{19}(m,p)=x_{\min}(p)+ v^{\rm meas}(p)\Big(-t_{\rm min}(p) +\\ \frac{x_{\min}(m)-x_{\min}(p)+v^{\rm meas}(p)t_{\min}(p)-vt_{\min}(m)}{v^{\rm meas}(p)-v} \Big)\\
t_{20}(m,p)=\frac{x_{\max}(m)-x_{\min}(p)+v^{\rm meas}(p)t_{\min}(p)-vt_{\max}(m)}{v^{\rm meas}(p)-w}\\
x_{20}(m,p)=x_{\min}(p)+v^{\rm meas}(p)\Big(-t_{\rm min}(p) +\\ \frac{x_{\max}(m)-x_{\min}(p)+v^{\rm meas}(p)t_{\min}(p)-vt_{\max}(m)}{v^{\rm meas}(p)-w} \Big)\\
t_{21}(m,p)=\frac{x_{\min}(m)-x_{\min}(p)+v^{meas}(p)t_{\min}(p)-vt_{\min}(m)}{v^{\rm meas}(p)-w}\\
x_{21}(m,p)=x_{\min}(p)+v^{\rm meas}(p)\Big(-t_{\rm min}(p)+\\ \frac{x_{\min}(m)-x_{\min}(p)+v^{meas}(p)t_{\min}(p)-vt_{\min}(m)}{v^{\rm meas}(p)-w} \Big)\\
x_{22}(m,u)=x_{\min}(m)+w(t_{\rho}(u)-t_{\min}(m)) \\
x_{23}(m,u)=x_{\max}(m)+w(t_{\rho}(u)-t_{\max}(m)) \\
x_{24}(m,u)=x_{\min}(m)+v(t_{\rho}(u)-t_{\min}(m)) \\
x_{25}(m,u)=x_{\max}(m)+v(t_{\rho}(u)-t_{\max}(m)) \\
x_{26}(m,u)=x_{\min}(m)+v_{meas}(m)(t_{\rho}(u)-t_{\min}(m)) \\
\end{array} \right.
\end{equation}

\begin{equation}\label{e:constants4}
\footnotesize
\left\{ \begin{array}{ll}
t_{27}(u,p)=\frac{x_{\min}(p)-x_{\max_{\rho}}(u)-v^{\rm meas}(p)t_{\min}(p)+vt_{\rho}(u)}{v-v^{\rm meas}(p)}\\
x_{27}(u,p)=x_{\min}(p)+v^{\rm meas}(p)\Big(-t_{\rm min}(p) +  \\  \frac{x_{\min}(p)-x_{\max_{\rho}}(u)-v^{\rm meas}(p)t_{\min}(p)+vt_{\rho}(u)}{v-v^{\rm meas}(p)}\Big)\\
t_{28}(u,p)=\frac{x_{\min}(p)-x_{\min_{\rho}}(u)-v^{\rm meas}(p)t_{\min}(p)+wt_{\rho}(u)}{w-v^{\rm meas}(p)}\\
x_{28}(u,p)=x_{\min}(p)+v^{\rm meas}(p)\Big(-t_{\rm min}(p)+ \\ \frac{x_{\min}(p)-x_{\min_{\rho}}(u)-v^{\rm meas}(p)t_{\min}(p)+wt_{\rho}(u)}{w-v^{\rm meas}(p)}  \Big)\\
t_{29}(u,p)=\frac{x_{\min}(p)-x_{\min_{\rho}}(u)-v^{\rm meas}(p)t_{\min}(p)+vt_{\rho}(u)}{v-v^{\rm meas}(p)}\\
x_{29}(u,p)=x_{\min}(p)+v^{\rm meas}(p)\Big(-t_{\rm min}(p)+ \\ \frac{x_{\min}(p)-x_{\min_{\rho}}(u)-v^{\rm meas}(p)t_{\min}(p)+vt_{\rho}(u)}{v-v^{\rm meas}(p)}  \Big)\\
t_{30}(u,p)=\frac{x_{\min}(p)-x_{\max_{\rho}}(u)-v^{\rm meas}(p)t_{\min}(p)+wt_{\rho}(u)}{w-v^{\rm meas}(p)}\\
x_{30}(u,p)=x_{\min}(p)+v^{\rm meas}(p)\Big(-t_{\rm min}(p)+ \\ \frac{x_{\min}(p)-x_{\max_{\rho}}(u)-v^{\rm meas}(p)t_{\min}(p)+wt_{\rho}(u)}{w-v^{\rm meas}(p)}  \Big)\\
x_{31}(u,p)=x_{\min_{\rho}}(u)+v(t_{\rho}(p)-t_{\rho}(u)) \\
x_{32}(u,p)=x_{\max_{\rho}}(u)+v(t_{\rho}(p)-t_{\rho}(u)) \\
x_{33}(u,p)=x_{\min_{\rho}}(u)+w(t_{\rho}(p)-t_{\rho}(u)) \\
x_{34}(u,p)=x_{\max_{\rho}}(u)+w(t_{\rho}(p)-t_{\rho}(u)) \\
\end{array} \right.
\end{equation}

{\bf Proof} --- \hspace{ 2 mm} Note that $\forall (k,n) \in [0, k_{\max}] \times [0,n_{\max}], {\rm Dom}(M_k) \cap {\rm Dom}({\bf M}_{\gamma_n})=\emptyset$ and that $\forall (k,n) \in [0, k_{\max}] \times [0,n_{\max}], {\rm Dom}(M_k) \cap {\rm Dom}({\bf M}_{\beta_n})=\emptyset$. Thus, the set of inequality constraints~(\ref{e:compconst}) can be written in the case of initial, upstream, downstream and internal conditions as:

\begin{equation} \label{e:ineqconsts1atemp}
\footnotesize
\begin{cases}
{\bf M}_{M_k}(t,\xi)\geq \gamma_p(t,\xi) & \forall t \in [pT,(p+1)T], \forall k\in\mathbb{K}, \forall p\in\mathbb{N}\hfill \\
{\bf M}_{M_k}(t,\chi)\geq \beta_p(t,\chi) & \forall t \in [pT,(p+1)T], \forall k\in\mathbb{K}, \forall p\in\mathbb{N}\hfill \\
{\bf M}_{M_k}(t,x) \geq \mu_m(t,x) & \forall (t,x) \in {\rm Dom}(\mu_m), \forall (k,m) \in \mathbb{K} \times \mathbb{M} \hfill \\
{\bf M}_{M_k}(t,x) \geq \Upsilon_u(t,x) & \forall (t,x) \in {\rm Dom}(\Upsilon_u), \forall (k,u) \in \mathbb{K} \times \mathbb{U} \hfill \\
{\bf M}_{\gamma_n}(t,\xi)\geq\gamma_p(t,\xi)  & \forall t \in [pT,(p+1)T], \forall (n,p)\in\mathbb{N}^2\hfill \\
{\bf M}_{\gamma_n}(t,\chi)\geq\beta_p(t,\chi)  & \forall t \in [pT,(p+1)T], \forall (n,p)\in\mathbb{N}^2\hfill \\
{\bf M}_{\gamma_n}(t,x) \geq \mu_m(t,x) & \forall (t,x) \in {\rm Dom}(\mu_m), \forall (n,m) \in \mathbb{N} \times \mathbb{M} \hfill \\
{\bf M}_{\gamma_n}(t,x) \geq \Upsilon_u(t,x) & \forall (t,x) \in {\rm Dom}(\Upsilon_u), \forall (k,u) \in \mathbb{K} \times \mathbb{U} \hfill \\
{\bf M}_{\beta_n}(t,\xi)\geq\gamma_p(t,\xi)  & \forall t \in [pT,(p+1)T], \forall (n,p)\in\mathbb{N}^2\hfill \\
{\bf M}_{\beta_n}(t,\xi)\geq\beta_p(t,\chi)  & \forall t \in [pT,(p+1)T], \forall (n,p)\in\mathbb{N}^2\hfill \\
{\bf M}_{\beta_n}(t,x) \geq \mu_m(t,x) & \forall (t,x) \in {\rm Dom}(\mu_m), \forall (n,m) \in \mathbb{N} \times \mathbb{M} \hfill \\
{\bf M}_{\beta_n}(t,x) \geq \Upsilon_u(t,x) & \forall (t,x) \in {\rm Dom}(\Upsilon_u), \forall (n,u) \in \mathbb{N} \times \mathbb{U} \hfill \\
{\bf M}_{\mu_k}(t,\xi)\geq \gamma_p(t,\xi) & \forall t \in [pT,(p+1)T], \forall k\in\mathbb{M}, \forall p\in\mathbb{N}\hfill \\
{\bf M}_{\mu_k}(t,\chi)\geq \beta_p(t,\chi) & \forall t \in [pT,(p+1)T], \forall k\in\mathbb{M}, \forall p\in\mathbb{N}\hfill \\
{\bf M}_{\mu_k}(t,x) \geq \mu_m(t,x) & \forall (t,x) \in {\rm Dom}(\mu_m), \forall (k,m) \in \mathbb{M}^2 \hfill \\
{\bf M}_{\mu_k}(t,x) \geq \Upsilon_u(t,x) & \forall (t,x) \in {\rm Dom}(\Upsilon_u), \forall k\in\mathbb{M}, \forall u\in\mathbb{U}\hfill  \\
{\bf M}_{\Upsilon_k}(t,\xi)\geq \gamma_p(t,\xi) & \forall t \in [pT,(p+1)T], \forall k\in\mathbb{U}, \forall p\in\mathbb{N}\hfill \\
{\bf M}_{\Upsilon_k}(t,\chi)\geq \beta_p(t,\chi) & \forall t \in [pT,(p+1)T], \forall k\in\mathbb{U}, \forall p\in\mathbb{N}\hfill \\
{\bf M}_{\Upsilon_k}(t,x) \geq \mu_m(t,x) & \forall (t,x) \in {\rm Dom}(\mu_m), \forall k \in \mathbb{U}, \forall m\in\mathbb{M} \hfill \\
{\bf M}_{\Upsilon_k}(t,x) \geq \Upsilon_u(t,x) & \forall (t,x) \in {\rm Dom}(\Upsilon_u), \forall (k,u)\in\mathbb{U}^2 \hfill  \\
\end{cases}
\end{equation}

The conditions~(\ref{e:ineqconsts1atemp}) all involve checking that a function of $(t,x)$ is greater than another function of $(t,x)$ on a line segment of $\mathbb{R}_+ \times [\xi,\chi]$.  Yet, because of the affine structure of the initial and boundary conditions~(\ref{e:defvalcondini}),~(\ref{e:defvalcondu}) and~(\ref{e:defvalcondd}) as well as the piecewise affine structure of their solutions~(\ref{e:solutions}),~(\ref{e:solutions2}) and~(\ref{e:solutions3}), the inequalities of the form $\forall (t,x) \in {\rm Dom}({\bf c}_i), \;\; {\bf M}_{{\bf c}_j}(t,x) \ge {\bf c}_i(t,x)$ are equivalent to a finite number of inequalities of the form $\forall p \in \{0,\dots,p_{\max}\}, \;\; {\bf M}_{{\bf c}_j}(t_p,x_p) \ge {\bf c}_i(t_p,x_p)$. This arises from the fact that a piecewise affine function is positive on all points of a segment if and only if it is positive on each extremity of the segment, and on the finite number of points of the segment on which the function is not differentiable. In the present case, this property implies the equivalence of~(\ref{e:ineqconsts1atemp}) and of~(\ref{e:ineqconsts1a}), (\ref{e:ineqconsts1b}), (\ref{e:ineqconsts1c}), (\ref{e:ineqconsts2a}), (\ref{e:ineqconsts2b}), (\ref{e:ineqconsts2c}), (\ref{e:ineqconsts3a}), (\ref{e:ineqconsts3b}), (\ref{e:ineqconsts3c}), (\ref{e:ineqconsts4a}), (\ref{e:ineqconsts4b}), (\ref{e:ineqconsts4c}), (\ref{e:ineqconsts4d}), (\ref{e:ineqconsts5a}), (\ref{e:ineqconsts5b}) and~(\ref{e:ineqconsts5c}). The equality constraints (\ref{e:eqconsts1}) arise from continuity conditions of the Moskowitz function~\cite{Newell93}. \hfill $\;\; \blacksquare$

\begin{Proposition} \label{p:continuityconst} \textbf{[Continuity constraints]} Let a set of initial, boundary and internal conditions be defined as in~(\ref{e:defvalcondini}), and let the corresponding partial solutions be defined as ${\bf M}_{M_k}(\cdot,\cdot)$, ${\bf M}_{\gamma_n}(\cdot,\cdot)$, ${\bf M}_{\beta_n}(\cdot,\cdot)$ and  ${\bf M}_{\mu_m}(\cdot,\cdot)$. Let us also assume that the model constraints~(\ref{e:compconst}) are satisfied. Let ${\bf M}_p(\cdot,\cdot)$ be defined as ${\bf M}_p(\cdot,\cdot)=\min_{k,n,u,m | m \neq p}({\bf M}_{M_k}(\cdot,\cdot),{\bf M}_{\gamma_n}(\cdot,\cdot),{\bf M}_{\beta_n}(\cdot,\cdot),{\bf M}_{\Upsilon_u}(\cdot,\cdot),\\{\bf M}_{\mu_m}(\cdot,\cdot))$ and ${\bf M}_o(\cdot,\cdot)$ be defined as ${\bf M}_o(\cdot,\cdot)=\min_{k,n,m,u | u \neq o}({\bf M}_{M_k}(\cdot,\cdot),{\bf M}_{\gamma_n}(\cdot,\cdot),{\bf M}_{\beta_n}(\cdot,\cdot),{\bf M}_{\mu_m}(\cdot,\cdot)),\\{\bf M}_{\Upsilon_u}(\cdot,\cdot)$ . The solution ${\bf M}(\cdot,\cdot)$ to the HJ PDE~(\ref{e:definfmorp}) defined by ${\bf M}(\cdot,\cdot)=\min_{k,n,m,u}({\bf M}_{M_k}(\cdot,\cdot),{\bf M}_{\gamma_n}(\cdot,\cdot),{\bf M}_{\beta_n}(\cdot,\cdot),{\bf M}_{\mu_m}(\cdot,\cdot),\\{\bf M}_{\Upsilon_u}(\cdot,\cdot))$ is continuous if and only if the following conditions are satisfied:

\begin{equation}{} \label{e:continuityconst}
\forall p \in \mathbb{M}, \; {\bf M}_p(t_{\rm min}(p),x_{\rm min}(p))=\mu_p(t_{\rm min}(p),x_{\rm min}(p))
\end{equation}

\begin{equation} \label{e:eqconsts1}
\footnotesize
\begin{cases}
\mu_m(t_{\min}(m),x_{\min}(m)) = \min\big({\bf M}_{M_k}(t_{\min}(m),x_{\min}(m)), \\{\bf M}_{\gamma_n}(t_{\min}(m),x_{\min}(m)), {\bf M}_{\beta_n}(t_{\min}(m),x_{\min}(m)), \\{\bf M}_{\Upsilon_u}(t_{\min}(m),x_{\min}(m)),{\bf M}_{\mu_p}(t_{\min}(m),x_{\min}(m))\big)  \hfill \forall k\in \mathbb{K},\\ \forall n\in \mathbb{N}, \forall u\in \mathbb{U},  \forall(m,p)\in\mathbb{M}^{2}\quad \\
\mu_m(t_{\max}(m),x_{\max}(m)) = \min\big({\bf M}_{M_k}(t_{\max}(m),x_{\max}(m)), \\{\bf M}_{\gamma_n}(t_{\max}(m),x_{\max}(m)), {\bf M}_{\beta_n}(t_{\max}(m),x_{\max}(m)), \\{\bf M}_{\Upsilon_u}(t_{\max}(m),x_{\max}(m)),{\bf M}_{\mu_p}(t_{\max}(m),x_{\max}(m))\big)  \hfill \forall k\in \mathbb{K},\\ \forall n\in \mathbb{N}, \forall u\in \mathbb{U},  \forall(m,p)\in\mathbb{M}^{2}\quad \\
\end{cases}
\end{equation}

\begin{equation}{} \label{e:continuityconst2}
\forall p \in \mathbb{M}, \; {\bf M}_o(t_{\rho}(o),x_{\rm min_{\rho}}(o))=\Upsilon_o(t_{\rho}(o),x_{\rm min_{\rho}}(o))
\end{equation}

\begin{equation} \label{e:eqconsts2}
\footnotesize
\begin{cases}
\Upsilon_u(t_{\rho}(u),x_{\min_{\rho}}(u)) = \min\big({\bf M}_{M_k}(t_{\rho}(u),x_{\min_{\rho}}(u)), \\{\bf M}_{\gamma_n}(t_{\rho}(u),x_{\min_{\rho}}(u)), {\bf M}_{\beta_n}(t_{\rho}(u),x_{\min_{\rho}}(u)), \\{\bf M}_{\mu_p}(t_{\rho}(u),x_{\min_{\rho}}(u)),{\bf M}_{\Upsilon_o}(t_{\rho}(u),x_{\min_{\rho}}(u))\big)  \hfill \forall k\in \mathbb{K},\\ \forall n\in \mathbb{N}, \forall p\in \mathbb{M},  \forall(o,u)\in\mathbb{U}^{2}\quad \\
\Upsilon_u(t_{\rho}(u),x_{\max_{\rho}}(u)) = \min\big({\bf M}_{M_k}(t_{\rho}(u),x_{\max_{\rho}}(u)), \\{\bf M}_{\gamma_n}(t_{\rho}(u),x_{\max_{\rho}}(u)), {\bf M}_{\beta_n}(t_{\rho}(u),x_{\max_{\rho}}(u)), \\{\bf M}_{\mu_p}(t_{\rho}(u),x_{\max_{\rho}}(u)),{\bf M}_{\Upsilon_o}(t_{\rho}(u),x_{\max_{\rho}}(u))\big)  \hfill \forall k\in \mathbb{K},\\ \forall n\in \mathbb{N}, \forall p\in \mathbb{M},  \forall(o,u)\in\mathbb{U}^{2}\quad \\
\end{cases}
\end{equation}

Furthermore, the equality constraints~(\ref{e:continuityconst}) and~(\ref{e:continuityconst2})  can be written as a set of mixed integer linear inequalities involving the continuous variables $\rho(1),\rho(2),\dots,\rho(k_{\max})$, $q_{\rm in}(1),\dots,q_{\rm in}(n_{\max}),$ $q_{\rm out}(1),\dots,q_{\rm out}(n_{\max})$, $L_1,\dots,L_{m_{\max}}$ and $r_1,\dots,r_{m_{\max}}$, as well as auxiliary integer variables.

\end{Proposition}

The proof of~(\ref{e:continuityconst}) is straightforward, and follows directly~\cite{SIAM11CB} from the piecewise affine structure of the partial solutions ${\bf M}_{M_k}(\cdot,\cdot)$, ${\bf M}_{\gamma_n}(\cdot,\cdot)$, ${\bf M}_{\beta_n}(\cdot,\cdot)$ and  ${\bf M}_{\mu_m}(\cdot,\cdot)$.

The fact that~(\ref{e:continuityconst}) can be written as a set of mixed integer linear inequalities is more involved. It can be shown that since ${\bf M}_p(\cdot,\cdot)=\min_{k,n,m|m \neq p}({\bf M}_{M_k}(\cdot,\cdot),$ ${\bf M}_{\gamma_n}(\cdot,\cdot),{\bf M}_{\beta_n}(\cdot,\cdot),{\bf M}_{\mu_m}(\cdot,\cdot))$, equation~(\ref{e:continuityconst}) can be written as a set of inequalities involving the continuous variables $\rho(1),\rho(2),\dots,\rho(k_{\max})$, $q_{\rm in}(1),\dots,q_{\rm in}(n_{\max}),$ $q_{\rm out}(1),\dots,q_{\rm out}(n_{\max})$, $L_1,\dots,L_{m_{\max}}$ and $r_1,\dots,r_{m_{\max}}$, as well as boolean variables. An example of such derivation is shown in~\cite{CCITS12} for the case in which $m_{\max}=1$. These inequalities can be further rewritten as mixed integer linear inequalities using the piecewise affine dependency of the partial solutions with respect to the variables $\rho(1),\rho(2),\dots,\rho(k_{\max})$, $q_{\rm in}(1),\dots,q_{\rm in}(n_{\max}),$ $q_{\rm out}(1),\dots,q_{\rm out}(n_{\max})$, $L_1,\dots,L_{m_{\max}}$ and $r_1,\dots,r_{m_{\max}}$.\\

\begin{Proposition} \label{p:dataconst} \textbf{[Data constraints]} In the remainder of our article, we assume that the data constraints are convex in the decision variable $v$.
\end{Proposition}

Different choices of error models yield convex data constraints, such as the two examples outlined below.

{\bf Example of convex data constraints (1)} --- \hspace{ 1 mm}  Consider a sensor measuring the boundary flows $(q_{\rm in}(0),...q_{\rm in}(n_{\max}))$ with $5\%$ relative uncertainty, a loop detector measuring the initial density $\rho(3)$ with $10 \%$ absolute uncertainty, and no downstream sensor. In this situation, the constraints are convex inequalities (indeed, linear inequalities) in the decision variable:

\begin{equation}
\footnotesize
\left\{ \begin{array}{cc}
0.95 q_{\rm in}^{\rm measured}(n) \le q_{\rm in}(n) \le 1.05 q_{\rm in}^{\rm measured}(n) \;\;\forall n \in [0,n_{\max}]\\
\rho(3)^{\rm measured} -0.1 \rho_m \le \rho(3) \le \rho(3)^{\rm measured} +0.1 \rho_m
\end{array} \right.
\end{equation}
\hfill $\;\; \blacksquare$

{\bf Example of convex data constraints (2)} --- \hspace{ 1 mm}  Consider two identical sensors measuring the boundary flows $(q_{\rm in}(0),...q_{\rm in}(n_{\max}))$, and $(q_{\rm out}(1),...q_{\rm out}(n_{\max}))$ which are characterized by a RMS relative error of $3\%$. In this situation, the constraints are convex inequalities (quadratic convex inequalities) in the decision variable:

\begin{equation}
\footnotesize
\left\{ \begin{array}{cc}
\displaystyle{\sum_{n=0}^{n_{\max}} \left( q_{\rm in}(n)-q_{\rm in}^{\rm measured} \right)^2} \le 0.03\displaystyle{\sum_{n=0}^{n_{\max}} \left(q_{\rm in}^{\rm measured} \right)^2 }\\
\displaystyle{\sum_{n=0}^{n_{\max}} \left( q_{\rm out}(n)-q_{\rm out}^{\rm measured} \right)^2} \le 0.03\displaystyle{\sum_{n=0}^{n_{\max}} \left(q_{\rm out}^{\rm measured} \right)^2 }
\end{array} \right.
\end{equation}

In this situation the estimation problem becomes a Quadratic Program.
\hfill $\;\; \blacksquare$

\begin{Proposition} \label{p:traveltimeconstdef} \textbf{[Travel time constraints]} Travel time data can be used in the estimation process, in order to properly define this information we define the travel time as:

\begin{equation}
\footnotesize
T_{time} = t_{f_{travel}} - t_{0_{travel}}
\end{equation}

The travel time constraints are specified as the following equality:

\begin{equation}
\footnotesize
\begin {cases}
{\bf M}_{\gamma_n}(t_{0_{travel}},\xi) = {\bf M}_{\beta_p}(t_{f_{travel}},\chi) &  \forall (n,p)\in\mathbb{N}^2\\
\end {cases}
\end{equation}
\end{Proposition}

\section{Traffic Estimation on a single highway link}
\label{s:trafficestsingle}

We now present an implementation of the estimation framework presented earlier on an experimental dataset. The dataset includes fixed sensor data (obtained from inductive loop detectors in the present case) travel time and mobile sensor data.

\subsection{Experimental setup}

In the following sections, the effectiveness of the method is illustrated on different traffic flow estimation problems which are formulated as LPs and MIPs, using the \emph{Mobile Century}~\cite{MMhhtp,ARMXWB09} dataset. The \emph{Mobile Century} field experiment demonstrated the use of \texttt{Nokia} N-95 cellphones as mobile traffic sensors in $2008$, and was a joint \texttt{UC-Berkeley/Nokia} project. \\

\begin{figure}[h]
\begin{center} \begin{tabular}{c}
{\mbox{\includegraphics[width = 0.9\linewidth]{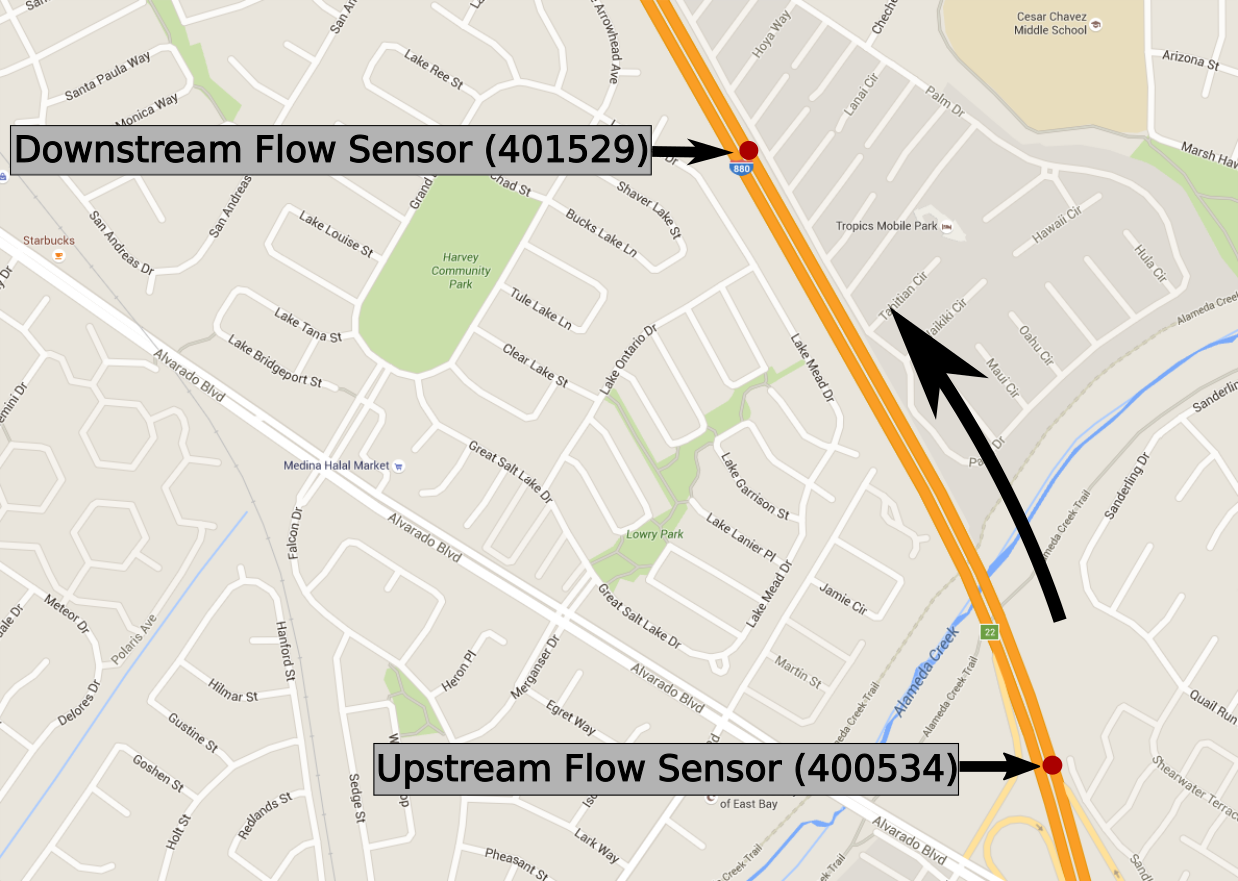}}}
\end{tabular}
\end{center}
\vspace{-0.1in} \caption{\textbf{Spatial domain considered for the numerical implementation.} The upstream and downstream PeMS stations are delimiting a $1.2\;km$ spatial domain, outlined by a solid line. The direction of traffic flow is represented by an arrow.\label{f:illlayout}}
\end{figure}

For the numerical applications, a spatial domain of 1.2 km is considered, located between the PeMS~\cite{PeMShttp} VDS stations 400536 and 401529 on the Highway I - 880 N around Hayward, California. The data used in this implementation was generated on February $8^{th}$, 2008, between times $18:30$ and $18:55$. In our scenario, we only consider inflow and outflow data $q_{in}^{\rm measured}(\cdot)$ and $q_{out}^{\rm measured}(\cdot)$ generated by the above PeMS stations, \emph{i.e.} we do not assume to know any density data. Of course the framework presented in this article allows incorporation of density data, see for instance Example 1 in the previous section. The layout of the spatial domain is illustrated in Figure~\ref{f:illlayout}.\\

\subsection{Initial density estimation on systems modeled by the Lighthill-Whitham-Richards PDE}

In this first scenario, our objective is to find the minimal or maximal values of a function of the decision variables, assuming that boundary flow data is available from the PeMS sensors. For this specific application the objective function is chosen as the total number of vehicles at initial time, defined by $\sum_{i=0}^{k_{max}}{\rho(i)}$, though any convex piecewise affine function of the decision variable would be acceptable. The constraints are linear inequalities in ~(\ref{e:defdecv}), and comprise both model~(\ref{e:ineqconsts1a}),~(\ref{e:ineqconsts2a}) and~(\ref{e:ineqconsts2c}) as well as data constraints. For this specific application, the data constraints are $(1-e)q_{\rm in/out}^{\rm measured}(n) \le q_{\rm in/out}(n) \le (1+e)q_{\rm in/out}^{\rm measured}(n) \;\;\forall n \in [0,n_{\max}]$, where $e=0.01=1\%$ is chosen the worst-case relative error of the sensors. The maximal densities are solutions to the following linear program:

\begin{equation} \label{e:datamodelcompat}
\footnotesize
\begin{array}{ll}
 \textbf{Minimize}\;&(\rm{respectively} \;\textbf{Maximize})\;  \sum_{i=0}^{k_{max}}{\rho(i)} \\
\\
\textbf{such that} & \left\{ \begin{array}{ll}(\ref{e:ineqconsts1a})\\(\ref{e:ineqconsts2a})\\(\ref{e:ineqconsts2c}) \\
(1-e)q_{\rm in}^{\rm measured}(n) \le q_{\rm in}(n) \;\;\forall n \in [0,n_{\max}]\\ q_{\rm in}(n) \le (1+e)q_{\rm in}^{\rm measured}(n) \;\;\forall n \in [0,n_{\max}]\\
(1-e)q_{\rm out}^{\rm measured}(n) \le q_{\rm out}(n) \;\;\forall n \in [0,n_{\max}]\\ q_{\rm out}(n) \le (1+e)q_{\rm out}^{\rm measured}(n) \;\;\forall n \in [0,n_{\max}]\\
\end{array} \right.
\end{array}
\end{equation}

\begin{figure}[h]
\begin{center} \begin{tabular}{c}
{\mbox{\includegraphics[width = 1.0\linewidth]{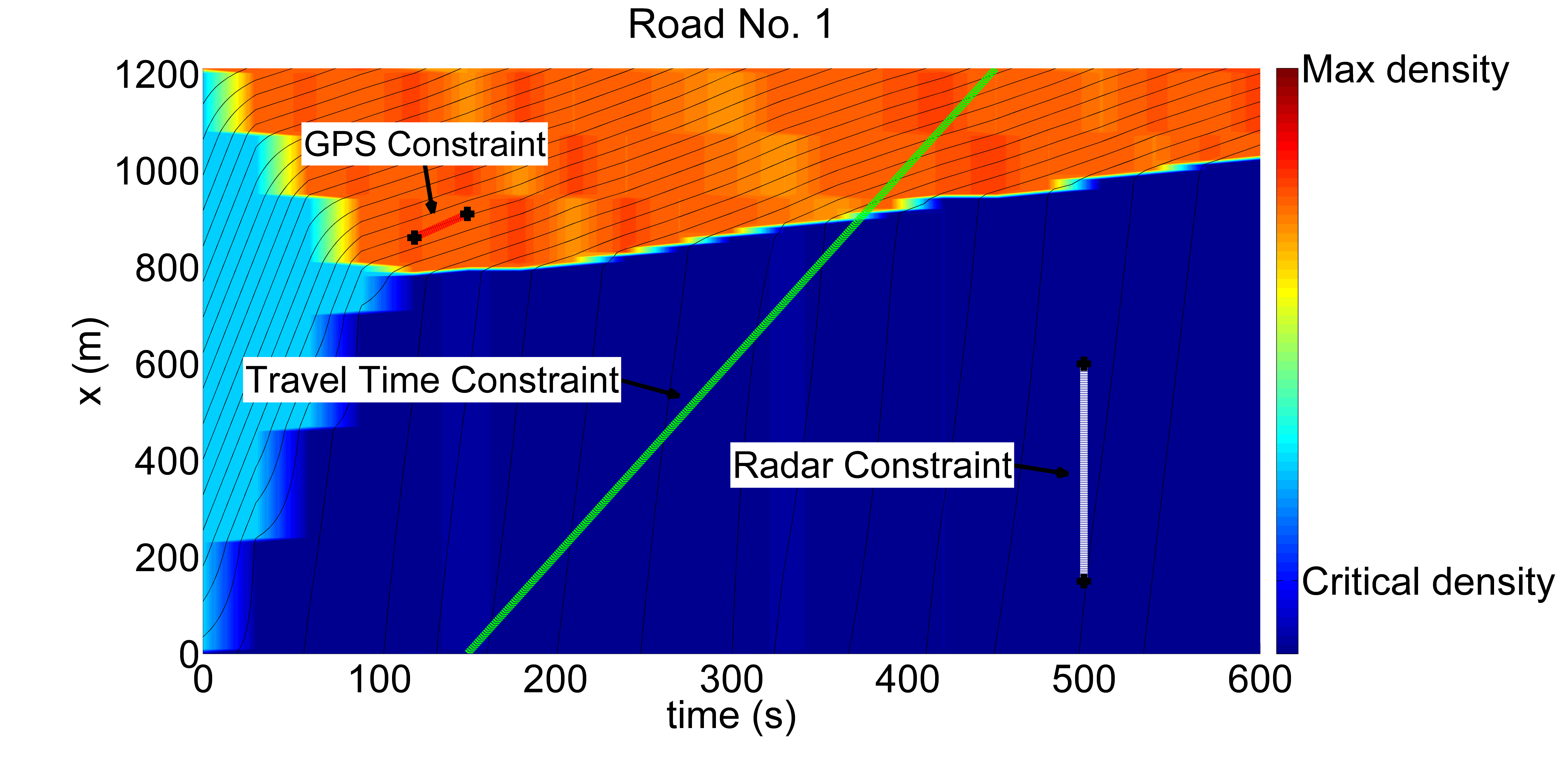}}}
\end{tabular}
\end{center}
\vspace{-0.1in} \caption{\textbf{Single road traffic state estimation.} In this problem, we want to reconstruct the state of traffic using three different types of data: GPS data generated by a probe vehicle, a travel time measurement generated by a vehicle entering and exiting the physical domain, and a radar generating a density measurement.  The present objective was to maximize the initial average density (worst-case average density), and the problem involves 49 variables and 929 constraints.\label{f:single_road}}
\end{figure}

For this implementation, we choose $20$ pieces of upstream and downstream flow data corresponding to $10$ minutes of data. The parameters of the triangular flux function are set with standard values: $v=65mph$, $w=-10mph$, $\rho_c=30 veh/(lane.mile)$.
The optimal solutions to~(\ref{e:datamodelcompat}) are illustrated in Figure~\ref{f:single_road}.

\section{Extension to Highway networks }
\label{s:highwaynet}

\subsection{Junction models}
In this section, we generalize the above framework  to traffic state estimation on road networks. For this, we first need to derive the boundary flows occurring at the junctions. This process is done through a junction model, which we now outline.

\begin{Proposition} \label{p:networkdefinition} \textbf{[Conservation of vehicles]} We consider a general junction (without storage capacity) integrating on-ramps, off-ramps and incoming as well as outgoing links, and illustrated in Figure \ref{f:networknode}.  

\begin{figure}[h]
\begin{center} \begin{tabular}{c}
{\mbox{\includegraphics[width = 0.7\linewidth]{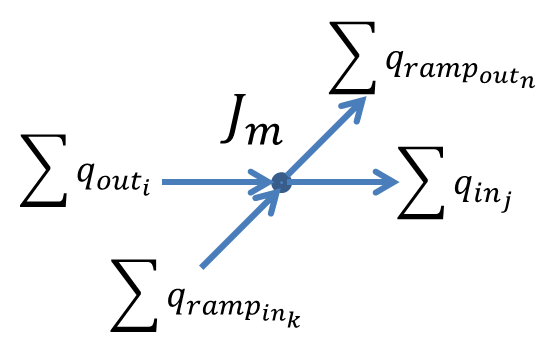}}}
\end{tabular}
\end{center}
\vspace{-0.1in} \caption{\textbf{Junction convention and flow conservation.}\label{f:networknode}}
\end{figure}

The conservation of vehicles at the junction imposes the following equality constraint.

\begin{equation} \label{e:conserv}
\footnotesize
\sum q_{out_{i}} + \sum q^{ramp}_{in_{k}} = \sum q_{in_{j}} + \sum q^{ramp}_{out_{n}}
\end{equation}

The quantities $\sum q_{out_{i}}$ are the flows leaving their respective road (to cross the junction) and $\sum q_{in_{j}}$ represent the flows entering the roads (after crossing the junction). $\sum q^{ramp}_{in_{k}}$ are the external flows entering the network through the intersection, conversely, $\sum q^{ramp}_{out_{n}}$ are the flows exiting the network.

\end{Proposition}

Note that the above equality constraint is insufficient in practice to derive the the actual incoming and outgoing flows. Following~\cite{coclite2005traffic,work2010traffic}, we assume that the flow at the junctions is the solution to a LP. This LP maximizes the total inflow through the junction under the constraints of an allocation matrix (which splits the flow according to the driver preferences), and physical constraints of demand and supply occuring at the boundaries of the incoming and outgoing links respectively.

The equation of conservation of flows~(\ref{e:conserv}) is linear in the decision variable. We assume that the volume of traffic entering the junction ($U_{i's}$) through $\sum q_{out_{i}}$ and $\sum q^{ramp}_{in_{k}}$ is distributed among the exit options ($O_{j's}$) according to an allocation parameter $\alpha_{Outs,Ins}\ge 0$. Since the junction has no storage capacity the sum of all the allocation parameters from a fixed incoming option among all the output options must be equal to one:

\begin{equation} \label{allocation}
\sum_{j} \alpha_{O_{j},U_{i}} = 1
\end{equation}

    The relation of the incoming-outgoing flows in terms of the allocation parameters is encoded by the allocation matrix:

\begin{equation} \label{flow}
\left[\sum_{j} O_{j} \right]=A \left[\sum_{i} U_{i}\right]
\end{equation}

The well-posedness constraints imposed by the LWR model on the junctions can be written as

\begin{equation} \label{e:demand}
\left\{ \begin{array}{l}
q_{out_{i}} \le d_{i}\\
 q^{ramp}_{in_{k}} \le d_{u,k}\\
q_{in_{j}} \le s_j\\
q^{ramp}_{out_{n}} \le s_{d,n}\\
\end{array} \right.
\end{equation}

where $d_{i}$ correspond to the demand of link $i$,  $d_{u,k}$ corresponds to the demand of the on-ramp $k$, $s_j$ correspond to the supply of link $j$ and 
$s_{d,n}$ correspond to the supply of off-ramp $n$. Since the demand and supplies of incoming and outgoing links is a mixed integer linear function of the decision variable~\cite{han2012link} (while the flows are linear in the decision variable), these constraints are mixed integer linear. Finally, the maximization of the sum of the flows through the junction imposes that at least one of the inequalities in~(\ref{e:demand}) becomes an equality, adding additional integer variables to the problem.

Hence, by combining all of the above constraints the junction constraints can be written as mixed integer linear inequalities in the decision variable, while the model constraints (on each links) are also mixed integer linear. If the objective to be minimized (as part of the estimation) is a linear function of the decision variable (which is typically the case), then the problem of estimating the state of traffic on a general highway network remains a MILP.

Hence, the estimation of the state of traffic on a general network can be posed as a mixed integer linear program, in which the integer variables are a consequence of the junction models and of internal data (either internal boundary conditions or internal density conditions).

\subsection{Implementation}

We illustrate the above results by implementing a MILP to solve two travel time estimation problems (with distinct network structures) involving different segments of the I-880 highway located in the San Francisco bay area. 

\begin{figure}[h]
\begin{center} \begin{tabular}{c}
{\mbox{\includegraphics[width = 1.0\linewidth]{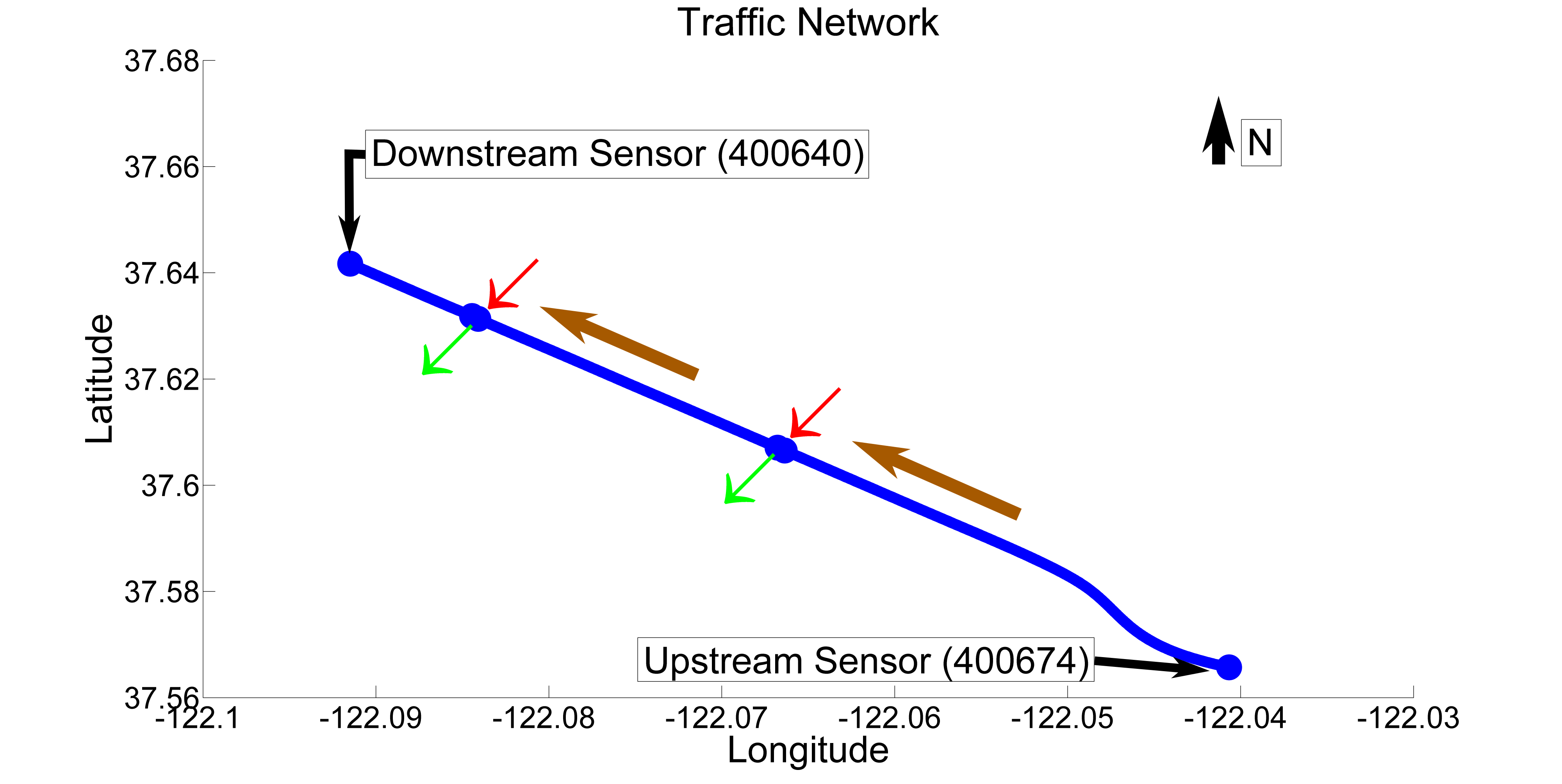}}}
\end{tabular}
\end{center}
\vspace{-0.1in} \caption{\textbf{Traffic Network Structure 1.} In this problem we want to estimate the travel time of the network from the upstream to downstream position; the network consists of three roads connected through two junctions with a ramp on and off respectively, the latter is used to consider the flow leaving and entering the higway on the intersections. \label{f:junction}}
\end{figure}

The first problem will involve three roads, two junctions with a ramp on and off respectively as can be seen in Figure \ref{f:junction}. For this problem, we consider data generated by the PeMS~\cite{PeMShttp} VDS stations 400674 and 400640 for the upstream and downstream position respectively (Figure \ref{f:junction}). The stations are located on the Highway I - 880 N around Hayward, California, the space domain of the entire network is 6 miles. We assumed that we don't have data in the junctions except for the allocation matrix which is predefined by the user, we are modeling the junctions with the approach described earlier. The traffic flow allocation matrix for both junctions in this example is the following:

\begin{equation} \label{matrix}
\left[ \begin{array}{l} 
q_{in_{1}} \\
 ramp_{out_{1}}  \end{array}\right] =\left[ \begin{array}{cc} 0.9 & 1  \\ 0.1 & 0 \end{array} \right] \begin{array}{l}  q_{out_{1}} \\ ramp_{in_{1}}  \end{array} 
\end{equation}

An example of the travel time estimation of the structure mentioned above is show in Figure~\ref{f:travel_time} below, the result obtained by the estimation toolbox is compared with the ground truth, that is, GPS data obtained during the \emph{Mobile Century} experiment. This particular example has 447 variables and 3127 constraints and took 4.57 seconds to solve on an iMac with a processor Intel Core i5-2400 @2.5GHz.

\begin{figure}[h]
\begin{center} \begin{tabular}{c}
{\mbox{\includegraphics[width = 1.0\linewidth]{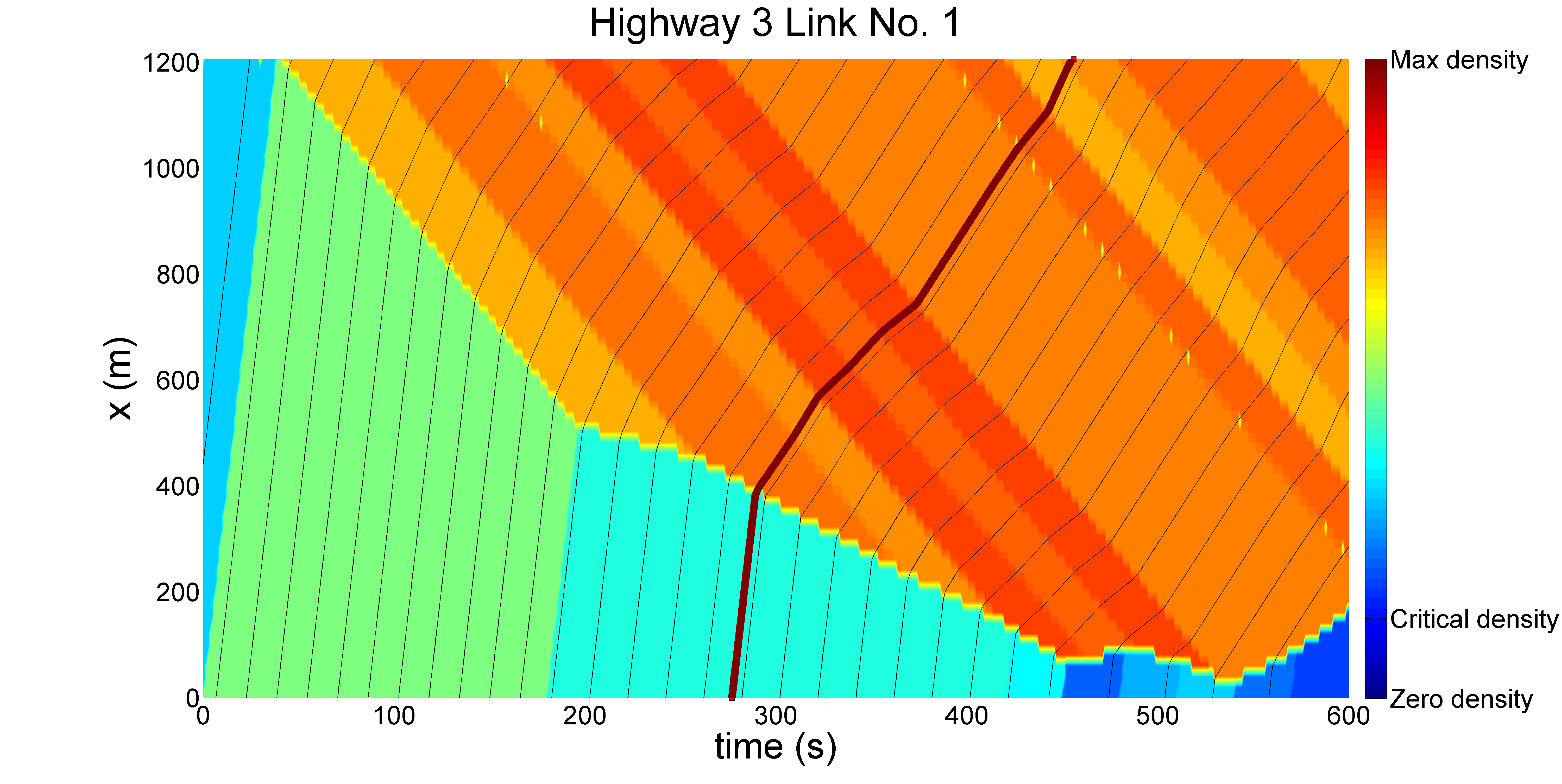}}}\\
{\mbox{\includegraphics[width = 1.0\linewidth]{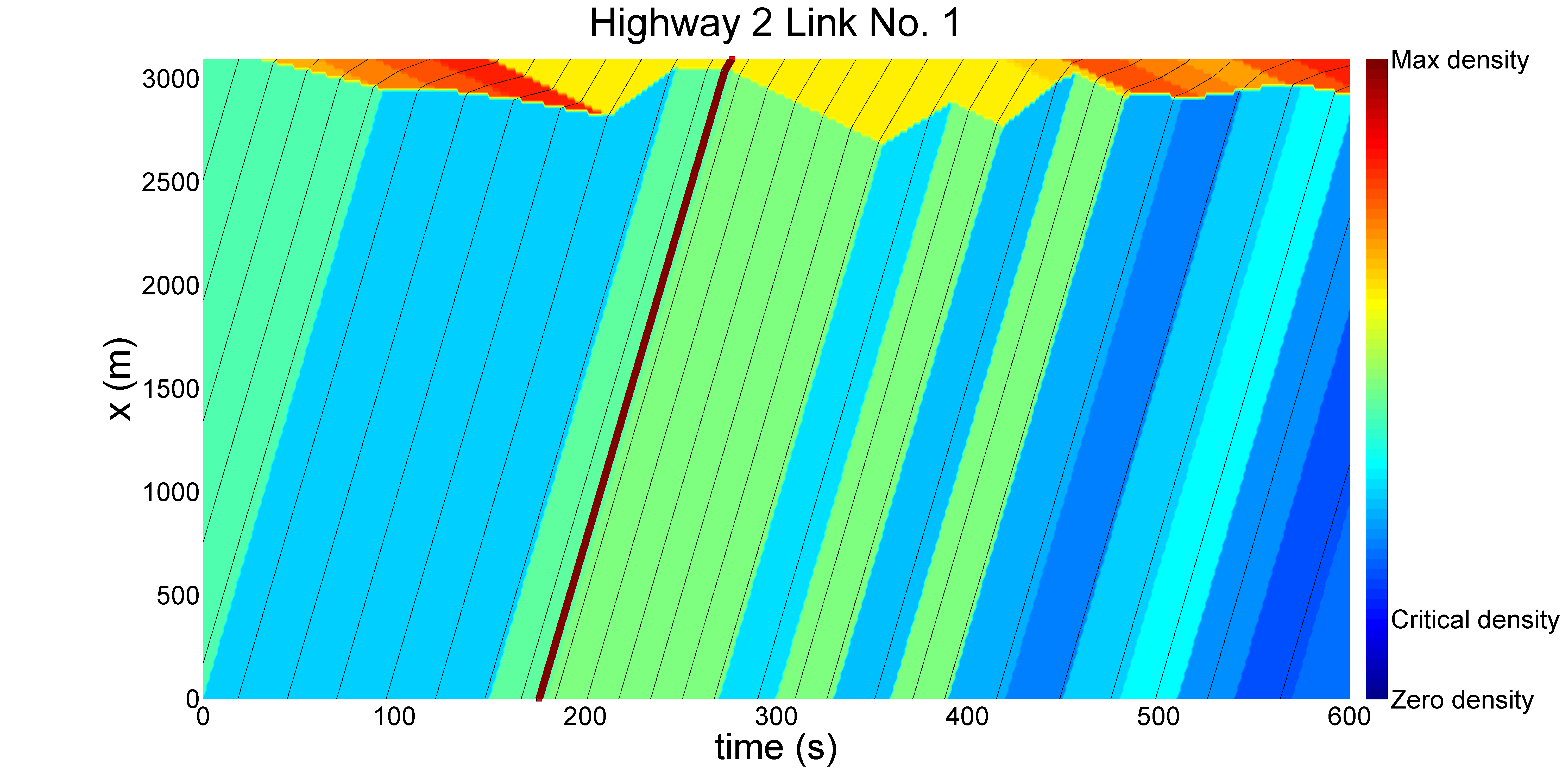}}}\\
{\mbox{\includegraphics[width = 1.0\linewidth]{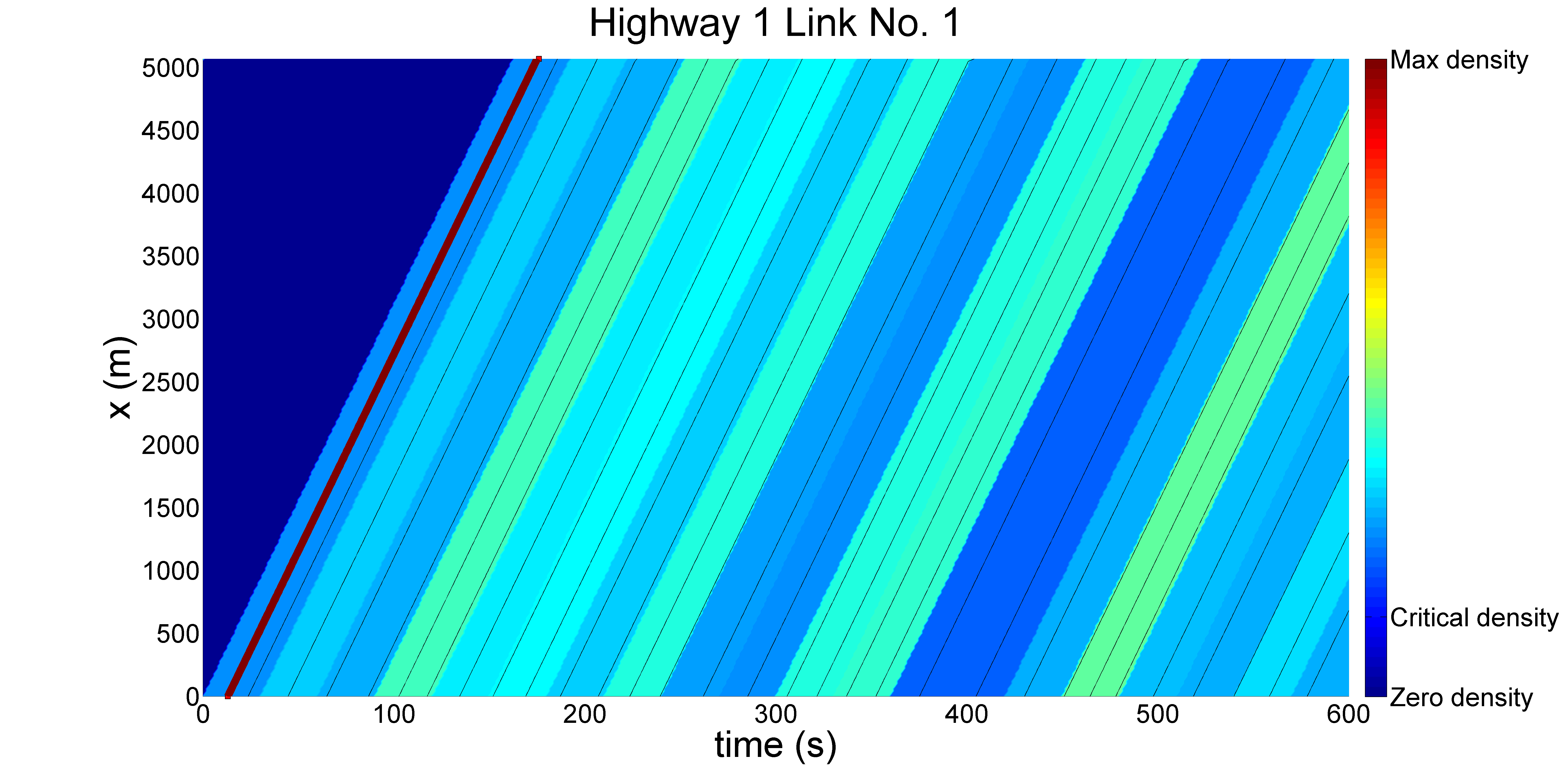}}}
\end{tabular}
\end{center}
\vspace{-0.1in} \caption{\textbf{Traffic Time Estimation Example Structure 1.} In all subfigures, we compute the density map for which the initial number of vehicles is the lowest, given the boundary data or the junction model, the sample car path is highlighted in red. For this example the estimated travel time is 443 seconds and the ground truth is 445 seconds. The objective function is the minimization of the initial densities. \textbf{Top:} Density map of the last highway segment, this segment has data in the downstream end. \textbf{Center:} Density map of the middle highway segment, this segment has a junction in both ends. \textbf{Bottom:} Density map of the initial segment of the highway, this segment has data on the upstream end. \label{f:travel_time}}
\end{figure}

A series of estimations with different traffic conditions were analyzed, comparing the estimated travel time with the ground truth obtained from GPS data. After 30 estimations the RMS error obtained is $11$ seconds. 

\begin{figure}[h]
\begin{center} \begin{tabular}{c}
{\mbox{\includegraphics[width = 1.0\linewidth]{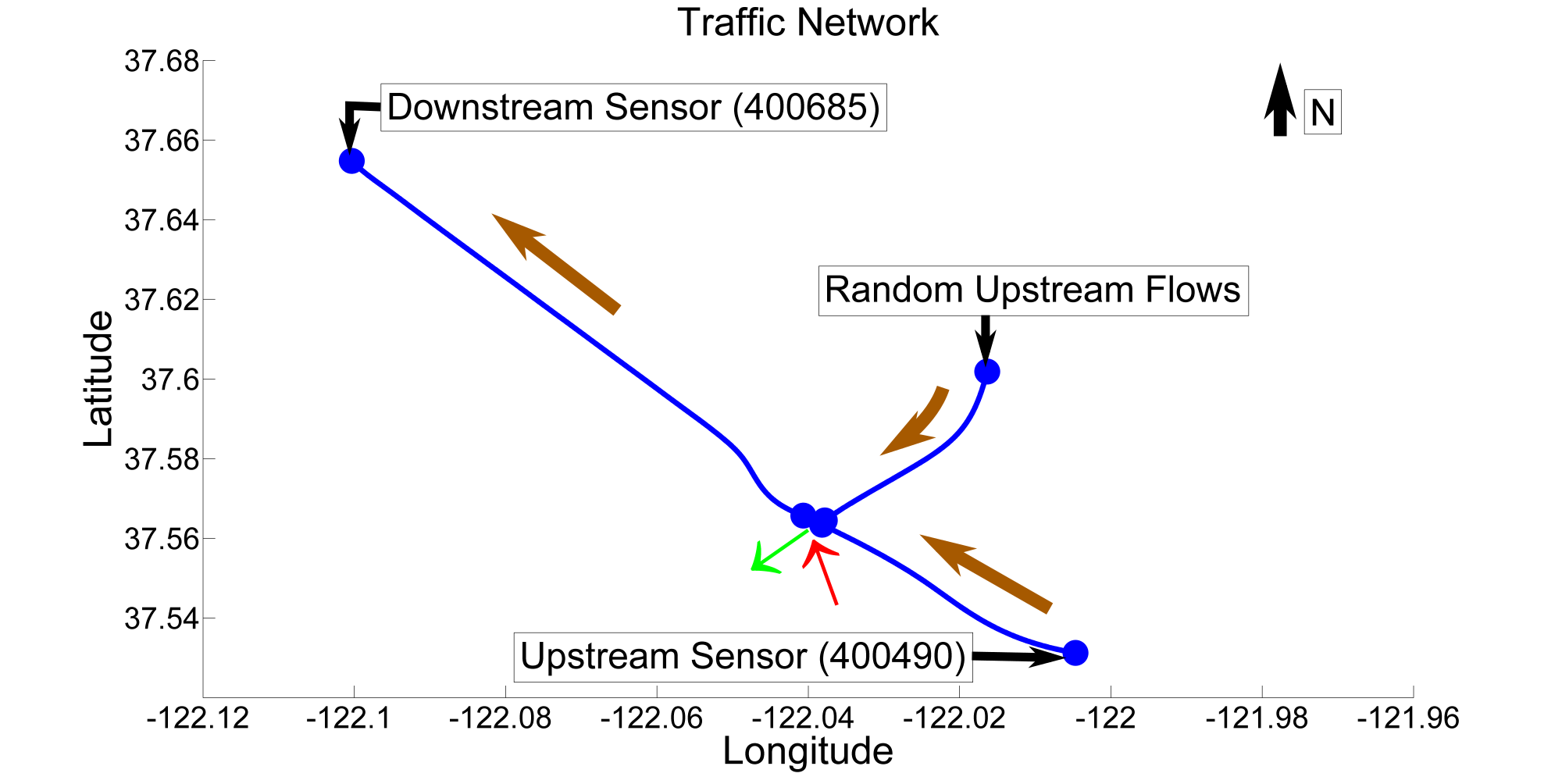}}}
\end{tabular}
\end{center}
\vspace{-0.1in} \caption{\textbf{Traffic Network Merge Structure.} We want to estimate the travel time of the network from the upstream to downstream position; the network consists of three roads connected through one junction with a ramp on and off, the latter is used to consider the flow leaving and entering the higway on the intersection. \label{f:junction2}}
\end{figure}

The second structure is a merge and will involve three roads and one junction with a ramp on and off as can be seen in Figure \ref{f:junction2}. For this problem, we consider data generated by the PeMS~\cite{PeMShttp} VDS stations 400490 and 400685 for the upstream and downstream position respectively (Figure \ref{f:junction}). The stations are located on the Highway I - 880 N around Hayward, California, the space domain of the entire network is 10.1 miles. Also, travel time data obtained on the \emph{Mobile Century} experiment was added as data constraints to the problem. The road merging to the I-880 is Decoto Road, from which there is no data available; a random generator of reasonable flows in the upstream position of Decoto Rd. was included. We assumed that we don't have data in the junction except for the allocation matrix which is again predefined by the user. The traffic flow allocation matrix in this problem is the following:

\begin{equation} \label{matrix}
\left[ \begin{array}{l} 
q_{in_{1}} \\
 ramp_{out_{1}}  \end{array}\right] =\left[ \begin{array}{ccc} 0.9 & 1 & 1 \\ 0.1 & 0 & 0 \end{array} \right] \begin{array}{l}  q_{out_{1}} \\ q_{out_{2}} \\ ramp_{in_{1}}  \end{array} 
\end{equation}

\begin{figure}[h]
\begin{center} \begin{tabular}{c}
{\mbox{\includegraphics[width = 1.0\linewidth]{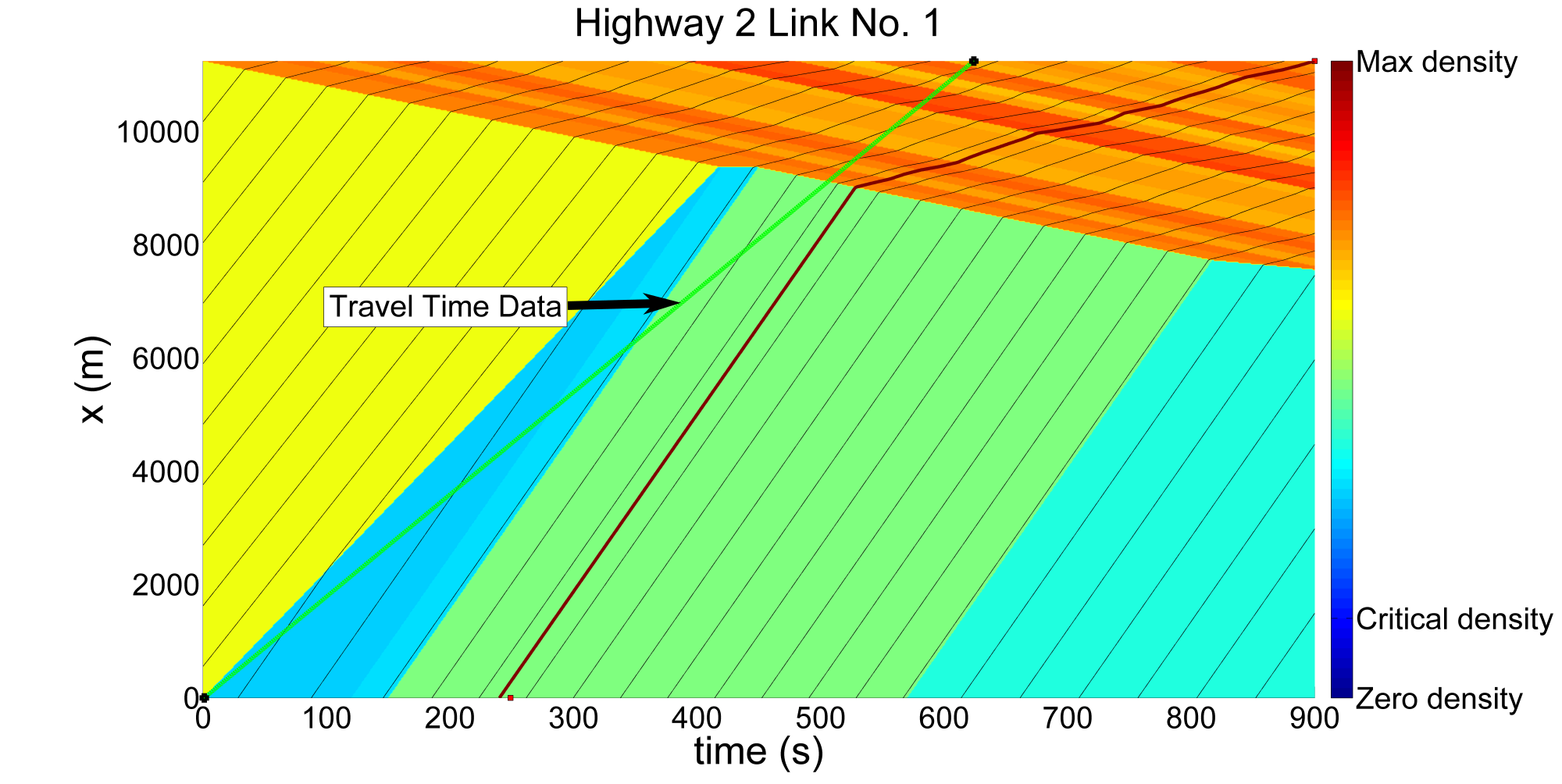}}}\\
{\mbox{\includegraphics[width = 1.0\linewidth]{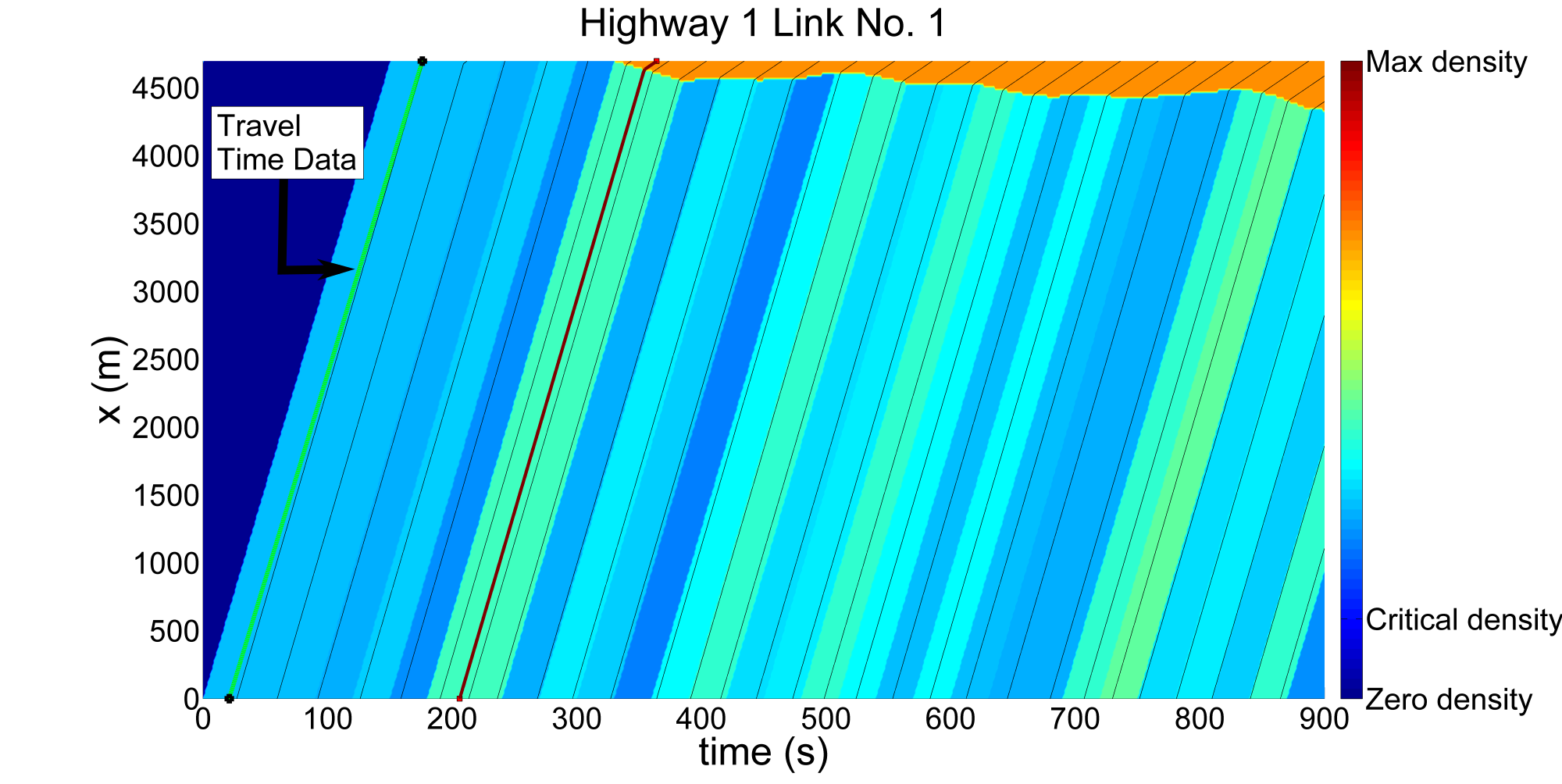}}}
\end{tabular}
\end{center}
\vspace{-0.1in} \caption{\textbf{Traffic Time Estimation Example 1 Merge Structure.} In these subfigures we compute the density map for which the initial number of vehicles is the lowest, given the boundary data or the junction model, the sample car path is highlighted in red. For this example the estimated travel time is 767 seconds and the ground truth is 769 seconds. \textbf{Top:} Density map of the last highway segment, this segment has data in the downstream end and a travel time constraint. \textbf{Bottom:} Density map of the initial segment of the highway, this segment has data on the upstream end and a travel time constraint. \label{f:travel_time_merge}}
\end{figure}

\begin{figure}[h]
\begin{center} \begin{tabular}{c}
{\mbox{\includegraphics[width = 1.0\linewidth]{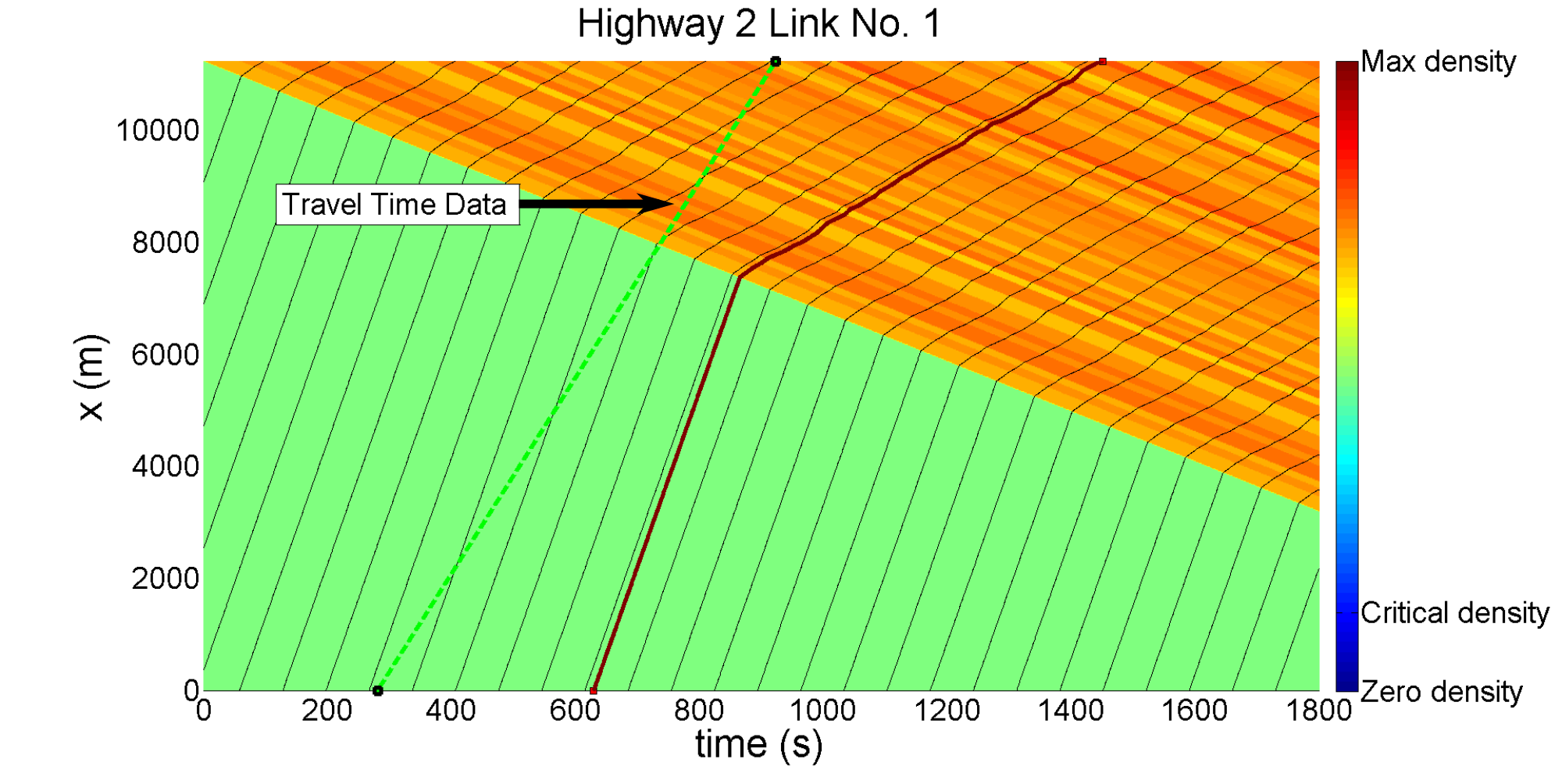}}}\\
{\mbox{\includegraphics[width = 1.0\linewidth]{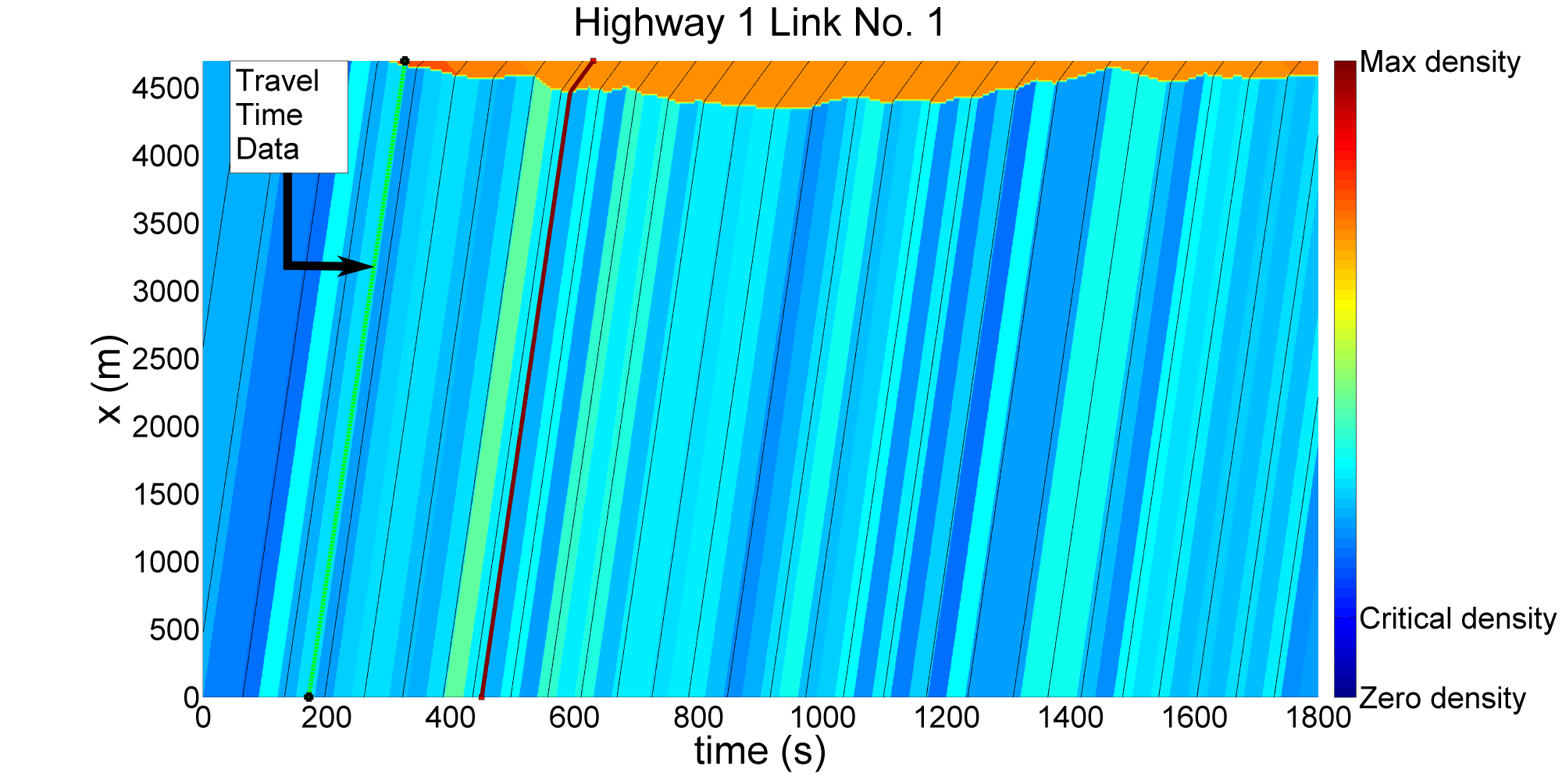}}}
\end{tabular}
\end{center}
\vspace{-0.1in} \caption{\textbf{Traffic Time Estimation Example 2 Merge Structure.} For these subfigures we compute the density map for which the L1 norm of the decision variables is the minimum, creating a more uniform density map. \textbf{Top:} Density map of the last highway segment, this segment has data in the downstream end and a travel time constraint. \textbf{Bottom:} Density map of the initial segment of the highway, this segment has data on the upstream end and a travel time constraint. \label{f:travel_time_merge2}}
\end{figure}

Two different objective functions were selected to estimate the travel time. The first objective function is the minimization of the initial densities, an example is shown in Figure~\ref{f:travel_time_merge}. The second objective function is the minimization of the L1 norm between the decision variables in order to obtain a more uniform density map (Figure \ref{f:travel_time_merge2}). The results obtained by the estimation toolbox are compared with the ground truth, that is, travel time data validated during the \emph{Mobile Century} experiment. The example inf Figure \ref{f:travel_time_merge} has 546 variables and 4336 constraints and took 5.12 seconds to solve on an iMac with a processor Intel Core i5-2400 @2.5GHz. The \texttt{MATLAB} toolbox can be downloaded from~\texttt{https://www.dropbox.com/s/} \texttt{1hyr9ffekb6731k/KAUST\_Traffic\_Network\_} \texttt{Estimation.zip?dl=0}.

\begin{figure}[h]
\begin{center} \begin{tabular}{c}
{\mbox{\includegraphics[width = 0.9\linewidth]{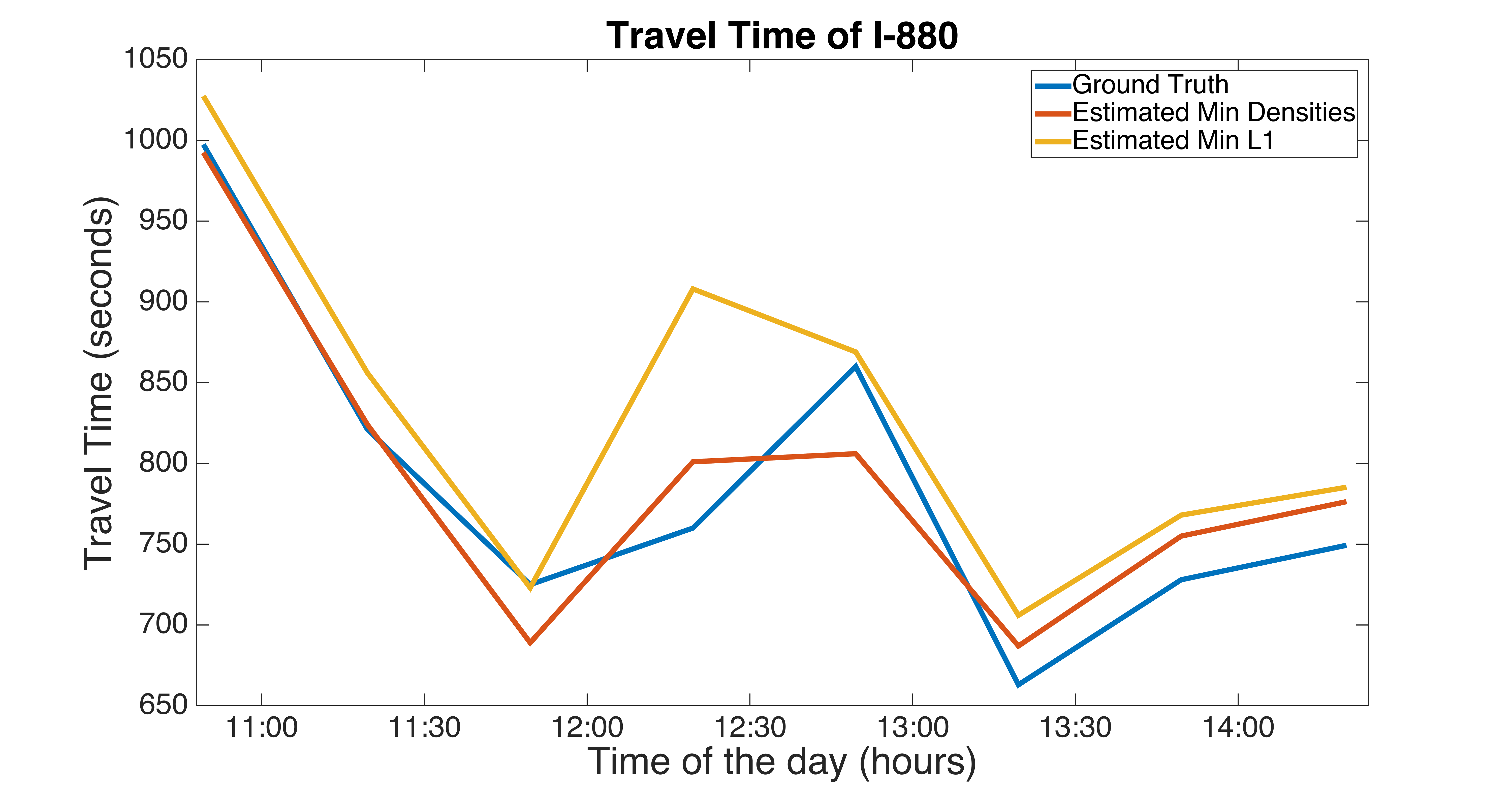}}}
\end{tabular}
\end{center}
\vspace{-0.1in} \caption{\textbf{Travel times Comparison.} The ground truth are the travel times experienced on February 8th, 2008 during the \emph{Mobile Century} experiment between Mowry and Winton Avenue (10.1 miles) compared with the results obtained by the travel time estimation using the merge structure and two different objective functions.\label{f:traveltimes}}
\end{figure}

The travel time estimation was evaluated using different traffic conditions, a comparison of the results obtained by the estimator with both objective functions can be observed in Figure~\ref{f:traveltimes}. From the figure we can conclude that with the L1 norm minimization the travel time is always overestimated, however it can capture better a sudden increase in the travel time, like the one happening at 13:00 hours of the evaluated day.

\section{Conclusion}

This article illustrates the travel time estimation application of a new Mixed Integer Programming estimation framework for highway traffic state estimation, in which the state of the system is modeled by the Lighthill Whitham Richards PDE. Using a Lax-Hopf formula, we show that the constraints arising from the model, as well as the the measurement data result in linear inequality constraints for a specific decision variable, and that the problem of estimating a linear function of this decision variable is a Mixed Integer Linear Program. We also show that the method can be generalized to highway networks, at the expense of increasing the decision variable size, hence affecting the computational time. A numerical implementation of the estimation on an experimental dataset containing fixed sensor data as well as probe data is performed, and illustrates the ability of the method to quickly and efficiently compute all traffic scenarios compatible both with both the LWR model and given traffic states. 

This framework has the advantage of being exact and efficient for small-scale networks and gives the ability for the user to select any objective function, to explore the possible state estimates associated with a given dataset. Future work will involve investigating a two phase flow model that has a realistic upper bound on acceleration, to model the state of traffic in urban environments. Preliminary analysis shows that this type of model could be integrated within a similar optimization framework.


\footnotesize
\bibliographystyle{plain}
\bibliography{biblioPDEcontrol,biblioHJE_LWR}

\end{document}